\documentclass{article}

\usepackage{arxiv}

\usepackage[utf8]{inputenc} 
\usepackage[T1]{fontenc}    
\usepackage{hyperref}       
\usepackage{url}            
\usepackage{booktabs}       
\usepackage{amsfonts}       
\usepackage{nicefrac}       
\usepackage{microtype}      
\usepackage{lipsum}		
\usepackage{graphicx}
\usepackage{natbib}
\usepackage{doi}
\usepackage{enumerate}
\usepackage{enumitem}
\usepackage{amsmath, amsthm, amssymb, algorithm, algpseudocode}
\usepackage[table]{xcolor}
\usepackage{tikz}

\newcommand{\gencontrol}{\mathfrak{M}}

\newtheorem{theorem}{Theorem}[section]
\newtheorem{lemma}[theorem]{Lemma}
\newtheorem{proposition}[theorem]{Proposition}
\newtheorem{corollary}[theorem]{Corollary}
\newtheorem{definition}[theorem]{Definition}
\newtheorem{example}[theorem]{Example}

\newtheorem{problem}[theorem]{Problem}

\newcommand\mystyle{\everymath{\displaystyle}}
\mystyle
\title{Second-Order $\Lambda$-Sets and Extensions to Non-Smooth, Hybrid, and Stochastic Optimal Control}

\author{\href{https://orcid.org/0000-0002-3816-5287}{\includegraphics[scale=0.06]{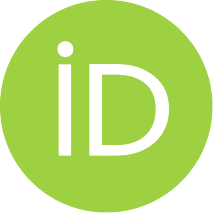}\hspace{1mm}M.H.M.~Rashid}\thanks{Corresponding Author} \\
	Department of Mathematics\&Statistics\\Faculty of Science P.O.Box(7)\\
	Mutah University University\\
	Mutah-Jordan \\
	\texttt{mrash@mutah.edu.jo}
}

\renewcommand{\shorttitle}{\textit{arXiv} Template}

\hypersetup{
pdftitle={Second-Order $\Lambda$-Sets and Extensions to Non-Smooth, Hybrid, and Stochastic Optimal Control},
pdfsubject={q-bio.NC, q-bio.QM},
pdfauthor={M.H.M.Rashid},
pdfkeywords={$\Lambda$-sets; optimal control; second-order conditions; Filippov systems; hybrid systems; stochastic hybrid systems; maximum principle; non-smooth dynamics}}

\begin{document}
\maketitle

\begin{abstract}
	 This paper develops a comprehensive extension of the $\Lambda$-set framework for optimal control, introducing second-order $\Lambda$-sets and generalizing the theory to non-smooth, hybrid, and stochastic hybrid systems. We first establish second-order necessary conditions that incorporate curvature information of the reachable set, providing refined optimality criteria that bridge classical second-variation methods with the geometric $\Lambda$-set approach. The framework is then extended to Filippov systems with discontinuous dynamics and to hybrid dynamical systems with state-dependent switching, yielding new necessary conditions for optimality in these settings. Furthermore, we introduce stochastic $\Lambda$-sets for systems subject to both continuous diffusion and discrete random switching, connecting the framework to Peng's stochastic maximum principle. Throughout the paper, detailed examples---including nonholonomic systems, mechanical systems with friction, and stochastic temperature control---illustrate the theoretical developments and demonstrate the practical applicability of the extended $\Lambda$-set theory. The results unify and generalize existing maximum principles, offering a powerful geometric tool for analyzing optimal control problems across a broad spectrum of system classes, from classical smooth systems to modern stochastic hybrid systems.
\end{abstract}

\section{Introduction}
\label{sec:introduction}

Optimal control theory represents a cornerstone of modern applied mathematics and engineering, providing powerful tools for optimizing dynamical systems across diverse domains including aerospace, robotics, economics, and biological systems. Since the pioneering work of Pontryagin and his collaborators \cite{Pontryagin1962}, the maximum principle has served as a fundamental necessary condition for optimality in deterministic control systems. However, classical formulations often rely on convexity assumptions and smooth dynamics, limiting their applicability to increasingly complex modern systems characterized by non-smooth behaviors, hybrid transitions, and stochastic disturbances.

The $\Lambda$-set framework, introduced and developed by Avakov and Magaril-Il'yaev \cite{Avakov2019, Avakov2020a, Avakov2020b}, represents a significant advancement in optimal control theory by providing a unified geometric approach to deriving necessary conditions for optimality without imposing traditional convexity requirements. This framework characterizes controllability properties through the existence of certain sets of adjoint variables, offering a bridge between geometric control theory and classical variational methods.

\subsection{Scholarly Contributions}

This paper makes several substantial contributions to the theory of $\Lambda$-sets and its applications:

\begin{enumerate}
    \item \textbf{Second-Order $\Lambda$-Sets:} We develop a comprehensive theory of second-order $\Lambda$-sets that incorporates curvature information of reachable sets. While first-order $\Lambda$-sets characterize local controllability through linear approximations, our second-order extension captures quadratic geometry, providing refined necessary conditions that bridge classical second-variation methods with the geometric $\Lambda$-set approach.

    \item \textbf{Extension to Non-Smooth Systems:} We extend the $\Lambda$-set framework to Filippov systems with discontinuous right-hand sides, establishing necessary conditions for optimality in systems with sliding modes and discontinuous dynamics. This extends the classical Filippov theory \cite{Filippov1959} to the $\Lambda$-set context.

    \item \textbf{Hybrid Systems Formulation:} We develop $\Lambda$-sets for hybrid dynamical systems with state-dependent switching, providing necessary conditions that unify continuous dynamics with discrete mode transitions. This contributes to the growing literature on hybrid optimal control \cite{Gamkrelidze1962, Lee1967}.

    \item \textbf{Stochastic Hybrid Systems:} We introduce stochastic $\Lambda$-sets for systems with both continuous diffusion and discrete random switching, connecting to Peng's stochastic maximum principle \cite{Peng1990} while extending it to hybrid settings.

    \item \textbf{Constructive Applications:} Through detailed examples including nonholonomic systems, mechanical systems with friction, and temperature control systems, we demonstrate both the theoretical validity and practical applicability of the extended $\Lambda$-set framework.
\end{enumerate}

\subsection{Significance and Motivation}

The significance of this research lies in several key aspects:

\subsubsection{Theoretical Significance}
\begin{itemize}
    \item \textbf{Unification of Frameworks:} The $\Lambda$-set approach provides a unified geometric perspective that encompasses classical maximum principles, second-order conditions, and extensions to non-smooth and hybrid systems. This unification offers deeper insight into the geometric structure of optimal control problems.

    \item \textbf{Gap Phenomenon Analysis:} The refined second-order conditions allow for precise characterization of the "gap phenomenon" \cite{FuscoMotta2022} where first-order conditions are satisfied but optimality fails, a situation commonly encountered in non-convex and non-smooth problems.

    \item \textbf{Connection to Controllability:} As established by Avakov and Magaril-Il'yaev \cite{Avakov2020b}, emptiness of $\Lambda$-sets implies local controllability, creating a direct link between optimality conditions and system controllability properties.
\end{itemize}

\subsubsection{Practical Motivation}
\begin{itemize}
    \item \textbf{Modern Engineering Applications:} Many contemporary engineering systems, from autonomous vehicles to robotic manipulators and power systems, exhibit hybrid behaviors, non-smooth dynamics, and stochastic disturbances. Traditional optimal control frameworks often struggle with these complexities, necessitating the extensions developed herein.

    \item \textbf{Computational Implications:} The $\Lambda$-set conditions provide verifiable certificates for optimality that can be incorporated into numerical algorithms for optimal control computation and verification.

    \item \textbf{Robustness Analysis:} The stochastic extensions provide tools for analyzing optimal control under uncertainty, a critical consideration in real-world applications where disturbances and modeling errors are inevitable.
\end{itemize}

\subsection{Applications and Relevance}

The extended $\Lambda$-set framework developed in this paper finds applications across multiple domains:

\subsubsection{Robotics and Autonomous Systems}
Nonholonomic systems, such as wheeled robots and autonomous vehicles, naturally give rise to control-affine systems where traditional convexity assumptions fail. The second-order $\Lambda$-set analysis enables rigorous treatment of time-optimal maneuvers and obstacle avoidance in such systems. The nonholonomic example in Section \ref{sec:second-order} demonstrates this application concretely.

\subsubsection{Mechanical Systems with Friction}
Systems with Coulomb friction and other discontinuous forces, as exemplified in Section \ref{sec:nonsmooth-hybrid}, require Filippov solutions and corresponding optimality conditions. Our extension of $\Lambda$-sets to Filippov systems provides necessary conditions for optimal control of mechanical systems with stick-slip phenomena and impact dynamics.

\subsubsection{Power Systems and Electrical Networks}
Hybrid systems naturally model power systems with switching components, circuit breakers, and mode changes. The hybrid $\Lambda$-set formulation enables optimal switching control for energy management, voltage regulation, and fault recovery.

\subsubsection{Biological and Chemical Systems}
Stochastic hybrid systems effectively model biochemical reactions, gene regulatory networks, and population dynamics where continuous dynamics interact with discrete events (e.g., cell division, chemical reactions). The stochastic hybrid $\Lambda$-sets developed in Section \ref{sec:stochastic} provide tools for optimal intervention in such systems.

\subsubsection{Temperature and Process Control}
The stochastic temperature control system analyzed in Section ~7 exemplifies industrial applications where optimal control must account for both stochastic disturbances (environmental fluctuations) and discrete mode changes (heating/cooling transitions).

\subsection{Relation to Existing Literature}

Our work builds upon and extends several important strands of research:

\begin{itemize}
    \item \textbf{Classical Maximum Principles:} The Pontryagin Maximum Principle \cite{Pontryagin1962} and its various extensions form the foundation upon which the $\Lambda$-set framework is constructed. Our second-order conditions relate to classical second-variation theory \cite{Zorich2004}.

    \item \textbf{$\Lambda$-Set Theory:} We extend the groundbreaking work of Avakov and Magaril-Il'yaev \cite{Avakov2019, Avakov2020a, Avakov2020b} to higher-order, non-smooth, and stochastic settings.

    \item \textbf{Non-Smooth Analysis:} Our treatment of Filippov systems connects to the extensive literature on differential inclusions and non-smooth analysis \cite{Filippov1959, Gamkrelidze1962}.

    \item \textbf{Hybrid Systems:} We draw upon hybrid optimal control theory \cite{Lee1967} and recent advances in impulsive control \cite{FuscoMotta2024, FuscoMotta2024LCSYS, FuscoMottaVinter2026}.

    \item \textbf{Stochastic Control:} The stochastic extensions connect to Peng's stochastic maximum principle \cite{Peng1990} and the theory of backward stochastic differential equations.
\end{itemize}

\subsection{Paper Organization}

The remainder of this paper is organized as follows: Section~\ref{sec:preliminaries} presents preliminary concepts, including the control system formulation, generalized controls, and the basic $\Lambda$-set framework; Section~\ref{sec:second-order} develops the theory of second-order $\Lambda$-sets, including second-order variations, quadratic expansions, and the main second-order necessary condition theorem; Section ~4 provides an example Nonholonomic System with Directional Constraint. Section~\ref{sec:nonsmooth-hybrid} extends the $\Lambda$-set framework to Filippov systems and hybrid dynamical systems, with applications to mechanical systems with friction and bouncing ball dynamics; Section~\ref{sec:stochastic} introduces stochastic $\Lambda$-sets for stochastic hybrid systems, establishing connections to classical stochastic maximum principles; Section~7 provides a detailed illustrative example of stochastic temperature control, demonstrating both theoretical and practical aspects of the framework; Section~8 identifies open problems and directions for future research; and the paper concludes with acknowledgments, declarations, and a comprehensive bibliography.

Through this systematic development, we aim to establish the extended $\Lambda$-set framework as a powerful and unified approach to optimal control across a broad spectrum of system classes, from classical smooth systems to modern stochastic hybrid systems.
\section{Preliminaries}
\label{sec:preliminaries}

This section establishes the foundational concepts, notation, and preliminary results that will be used throughout the paper. We begin by defining the class of control systems under consideration, followed by essential notions from convex analysis, generalized controls, and the basic $\Lambda$-set framework.

\subsection{Control System Formulation}

Consider a control system of the form:
\begin{equation}
\label{eq:control-system}
\dot{x}(t) = f(t, x(t), u(t)), \quad t \in [t_1, t_2],
\end{equation}
with initial condition
\begin{equation}
\label{eq:initial-condition}
x(t_1) = x_1,
\end{equation}
where:
\begin{itemize}
    \item $x(t) \in \mathbb{R}^n$ is the state at time $t$,
    \item $u(t) \in U \subset \mathbb{R}^r$ is the control input, with $U$ a compact set,
    \item $f: [t_1, t_2] \times \mathbb{R}^n \times U \to \mathbb{R}^n$ is a continuous mapping, Lipschitz in $x$ uniformly in $(t, u)$.
\end{itemize}

We denote by $\mathcal{U}$ the set of measurable control functions $u: [t_1, t_2] \to U$. For each $u \in \mathcal{U}$, there exists a unique absolutely continuous solution $x(\cdot; u)$ of \eqref{eq:control-system}--\eqref{eq:initial-condition} defined on $[t_1, t_2]$.

\subsection{Convexification and Generalized Controls}

The classical existence theory for optimal control problems often requires convexity assumptions on the velocity set $\{f(t, x, u): u \in U\}$. To handle non-convex problems, we employ the \emph{relaxed} or \emph{convexified} formulation using generalized controls \cite{Varga1962, Lee1967, Gamkrelidze1978}.

\begin{definition}[Generalized Control]
A \emph{generalized control} is a family $\mu = \{\mu_t\}_{t \in [t_1, t_2]}$ of Borel probability measures on $U$ such that for every Borel set $B \subset U$, the mapping $t \mapsto \mu_t(B)$ is measurable. We denote by $\mathfrak{M}_U$ the set of all generalized controls.
\end{definition}

The convexified dynamics corresponding to a generalized control $\mu \in \mathfrak{M}_U$ are given by:
\begin{equation}
\label{eq:convexified}
\dot{x}(t) = \int_U f(t, x(t), u) \, d\mu_t(u) = \langle \mu_t, f(t, x(t), \cdot) \rangle.
\end{equation}

For an ordinary control $u(\cdot) \in \mathcal{U}$, the corresponding generalized control is $\mu_t = \delta_{u(t)}$, the Dirac measure concentrated at $u(t)$. The set of trajectories generated by generalized controls contains the closure (in the uniform topology) of the set of trajectories generated by ordinary controls.

\subsection{Endpoint Mapping and Controllability}

Let $T > t_1$ and define the \emph{endpoint mapping} $\mathcal{E}_T: \mathfrak{M}_U \to \mathbb{R}^n$ by
\[
\mathcal{E}_T(\mu) = x(T; \mu),
\]
where $x(\cdot; \mu)$ solves \eqref{eq:convexified} with initial condition \eqref{eq:initial-condition}. The endpoint mapping is Fréchet differentiable under appropriate smoothness assumptions on $f$ \cite{Avakov2020b}.

\begin{definition}[Local Controllability]
The system \eqref{eq:control-system} is \emph{locally controllable} from $x_1$ at time $T$ if for every neighborhood $\mathcal{V}$ of the reference trajectory $\widehat{x}(\cdot)$ in the space of continuous functions, there exists $\tau \in (t_1, T)$ and an ordinary control $u \in \mathcal{U}$ such that the corresponding trajectory $x(\cdot; u) \in \mathcal{V}$ and $x(\tau; u) = \widehat{x}(T)$.
\end{definition}

\subsection{First-Order $\Lambda$-Sets}

The $\Lambda$-set framework, introduced by Avakov and Magaril-Il'yaev \cite{Avakov2019, Avakov2020a}, provides a geometric approach to deriving necessary conditions for optimality. Let $(\widehat{x}(\cdot), \widehat{\mu}_t)$ be an admissible pair for the convexified system.

\begin{definition}[First-Order $\Lambda$-Set $\Lambda_s$]
For $s \in \{-1, 1\}$, the \emph{first-order $\Lambda$-set} $\Lambda_s(\widehat{x}, \widehat{\mu})$ consists of absolutely continuous functions $\psi: [t_1, \widehat{t}_2] \to (\mathbb{R}^n)^*$ such that:
\begin{enumerate}
    \item $\psi$ satisfies the adjoint equation:
    \begin{equation}
    \label{eq:adjoint}
    \dot{\psi}(t) = -\psi(t) \langle \widehat{\mu}_t, f_x(t, \widehat{x}(t), u) \rangle \quad \text{a.e. } t \in [t_1, \widehat{t}_2].
    \end{equation}

    \item The maximum condition holds:
    \begin{equation}
    \label{eq:max-condition}
    \langle \psi(t), \dot{\widehat{x}}(t) \rangle = \max_{u \in U} H(t, \widehat{x}(t), \psi(t), u) \quad \text{a.e. } t,
    \end{equation}
    where $H(t, x, \psi, u) = \langle \psi, f(t, x, u) \rangle$ is the Hamiltonian.

    \item The transversality condition is satisfied:
    \begin{equation}
    \label{eq:transversality}
    s \, M(\widehat{t}_2, \widehat{x}(\widehat{t}_2), \psi(\widehat{t}_2)) \leq 0,
    \end{equation}
    with $M(t, x, \psi) = \max_{u \in U} H(t, x, \psi, u)$.

    \item $\psi \neq 0$ (nontriviality condition).
\end{enumerate}
\end{definition}

The parameter $s$ distinguishes between necessary conditions for minimality ($s = -1$) and maximality ($s = 1$). The following fundamental result links $\Lambda$-sets to optimality:

\begin{theorem}[First-Order Necessary Condition \cite{Avakov2020b}]
\label{thm:main}
If $\widehat{x}(\cdot)$ is a strong local minimum for the time-optimal problem, then for every corresponding generalized control $\widehat{\mu}_t$, the first-order $\Lambda$-set $\Lambda_{-1}(\widehat{x}, \widehat{\mu})$ is nonempty.
\end{theorem}

Conversely, emptiness of $\Lambda_{-1}$ implies that the reference trajectory cannot be a strong local minimum, and in fact, there exists a nearby trajectory reaching the target in shorter time.

\subsection{Notation and Basic Assumptions}

Throughout this paper, we adopt the following notation and assumptions:

\begin{itemize}
    \item $\mathbb{R}^n$ denotes $n$-dimensional Euclidean space with inner product $\langle \cdot, \cdot \rangle$ and norm $\|\cdot\|$.

    \item $\mathbb{S}^n$ denotes the space of symmetric $n \times n$ matrices.

    \item For a matrix $A$, $A^\top$ denotes its transpose, and $A \preceq B$ denotes the Loewner order ($B - A$ is positive semidefinite).

    \item $C^k$ denotes the space of $k$-times continuously differentiable functions, equipped with the standard norm.

    \item We assume $f \in C^2$ with respect to $x$ and $u$, with all derivatives bounded on compact sets.

    \item The control set $U$ is compact and convex (or has convex hull used in the convexified formulation).

    \item The reference trajectory $\widehat{x}(\cdot)$ is defined on $[t_1, \widehat{t}_2]$ and corresponds to a generalized control $\widehat{\mu}_t \in \mathfrak{M}_U$.
\end{itemize}

\subsection{Relationship to Classical Maximum Principles}

The $\Lambda$-set framework generalizes and refines classical necessary conditions for optimality. In particular:
\begin{itemize}
    \item When $f$ is smooth and the control set is convex, $\Lambda_{-1}$ nonemptiness implies the Pontryagin Maximum Principle \cite{Pontryagin1962}.

    \item For free end-time problems, the transversality condition \eqref{eq:transversality} incorporates the time optimization aspect, relating to the Hamiltonian value at the terminal time.

    \item The $\Lambda$-set approach naturally extends to non-smooth and hybrid systems, as will be shown in subsequent sections.
\end{itemize}

The following sections will build upon these preliminaries to develop second-order $\Lambda$-sets, extend the framework to non-smooth and hybrid systems, and finally to stochastic hybrid systems.
latex
\section{Second-Order $\Lambda$-Sets and Higher-Order Necessary Conditions}
\label{sec:second-order}

This section presents a novel extension of the $\Lambda$-set framework to second-order conditions,
establishing refined necessary conditions for optimality that incorporate curvature information
of the reachable set. While first-order $\Lambda$-sets characterize local controllability through
linear approximations, second-order $\Lambda$-sets capture the quadratic geometry of reachable sets
near optimal trajectories, providing a bridge between classical second-variation methods
and the geometric approach of $\Lambda$-set theory.

\subsection{Second-Order Variations and Quadratic Expansion}

Let $(\widehat{x}(\cdot), \widehat{\mu}_t)$ be an admissible pair for the convexified system
\begin{equation}
\label{eq:convex}
\dot{x} = \langle \mu_t, f(t, x, u) \rangle, \quad \mu_t \in \mathfrak{M}_U, \quad x(t_1) = x_1
\end{equation}
on $[t_1, \widehat{t}_2]$. Consider two-parameter variations of the form:
\[
\mu_t^{\alpha,\beta} = \widehat{\mu}_t + \alpha\delta\mu_t^1 + \beta\delta\mu_t^2 + \frac{\alpha^2}{2}\delta^2\mu_t^{11} + \alpha\beta\delta^2\mu_t^{12} + \frac{\beta^2}{2}\delta^2\mu_t^{22}
\]
where $\delta\mu_t^i \in \mathfrak{M}_U - \widehat{\mu}_t$ and $\delta^2\mu_t^{ij}$ are second-order
measure variations satisfying appropriate measurability conditions.
\begin{definition}[Second-Order Endpoint Mapping]
The \emph{second-order endpoint mapping} $\mathcal{E}^2: \mathbb{R}^2 \to \mathbb{R}^n$ is defined as:
\[
\mathcal{E}^2(\alpha,\beta) = x(\widehat{t}_2; \mu^{\alpha,\beta}) - \widehat{x}(\widehat{t}_2)
\]
where $x(\cdot; \mu^{\alpha,\beta})$ solves the convexified system with control $\mu^{\alpha,\beta}$.
\end{definition}

The Taylor expansion of $\mathcal{E}^2$ at $(0,0)$ yields:

\begin{align*}
\mathcal{E}^2(\alpha,\beta) &= \alpha \nabla_\alpha\mathcal{E}^2(0,0) + \beta \nabla_\beta\mathcal{E}^2(0,0) \\
&\quad + \frac{1}{2}\Bigl(\alpha^2 \nabla_{\alpha\alpha}^2\mathcal{E}^2(0,0) + 2\alpha\beta \nabla_{\alpha\beta}^2\mathcal{E}^2(0,0) \\
&\quad + \beta^2 \nabla_{\beta\beta}^2\mathcal{E}^2(0,0)\Bigr) + o(\alpha^2+\beta^2)
\end{align*}

The first-order terms correspond to the linearized system, while the second-order terms
involve solutions of the second variation equation.

\subsection{Second-Order Adjoint System and Curvature Form}

Consider the second variation system along $(\widehat{x},\widehat{\mu})$:
\begin{equation}
\begin{aligned}
\delta\dot{x} &= \langle \widehat{\mu}_t, f_x(t,\widehat{x},u)\delta x \rangle + \langle \delta\mu_t, f(t,\widehat{x},u) \rangle \\
\delta^2\dot{x} &= \langle \widehat{\mu}_t, f_x(t,\widehat{x},u)\delta^2 x \rangle + \langle \widehat{\mu}_t, f_{xx}(t,\widehat{x},u)(\delta x,\delta x) \rangle \\
&\quad + 2\langle \delta\mu_t, f_x(t,\widehat{x},u)\delta x \rangle + \langle \delta^2\mu_t, f(t,\widehat{x},u) \rangle
\end{aligned}
\end{equation}

\begin{definition}[Second-Order Hamiltonian]
For $(x,\psi,Q) \in \mathbb{R}^n \times (\mathbb{R}^n)^* \times \mathbb{S}^n$, where $\mathbb{S}^n$
denotes symmetric $n\times n$ matrices, define the \emph{second-order Hamiltonian}:
\begin{equation}\label{M1}
\begin{split}
H^2(t,x,\psi,Q,u) &= \langle \psi, f(t,x,u) \rangle \\
&\quad + \frac{1}{2}\langle Q f_x(t,x,u) + f_x(t,x,u)^\top Q, \delta x \otimes \delta x \rangle \\
&\quad + \langle Q, f_{xx}(t,x,u)(\delta x,\delta x) \rangle
\end{split}
\end{equation}

The \emph{second-order maximum function} is:
\begin{equation}
M^2(t,x,\psi,Q,\delta x) = \sup_{u \in U} H^2(t,x,\psi,Q,u,\delta x)
\end{equation}
\end{definition}

\subsection{Second-Order $\Lambda$-Set Definition}

\begin{definition}[Second-Order $\Lambda$-Set $\Lambda^{(2)}_s$]
For $s \in \{-1,1\}$, the \emph{second-order $\Lambda$-set} $\Lambda^{(2)}_s(\widehat{x},\widehat{\mu})$
consists of triples $(\psi, Q, \Psi)$ where:
\begin{enumerate}
\item $\psi: [t_1,\widehat{t}_2] \to (\mathbb{R}^n)^*$ is absolutely continuous and satisfies the
first-order adjoint equation:
\[
\dot{\psi} = -\psi \langle \widehat{\mu}_t, f_x(t,\widehat{x},u) \rangle
\]

\item $Q: [t_1,\widehat{t}_2] \to \mathbb{S}^n$ is absolutely continuous and satisfies the
matrix Riccati differential inequality:
\[
\dot{Q} + Q\langle \widehat{\mu}_t, f_x(t,\widehat{x},u) \rangle + \langle \widehat{\mu}_t, f_x(t,\widehat{x},u) \rangle^\top Q
+ \langle \widehat{\mu}_t, f_{xx}(t,\widehat{x},u) \rangle \preceq 0
\]
where $\preceq$ denotes the Loewner order on symmetric matrices.

\item $\Psi: [t_1,\widehat{t}_2] \to \mathbb{R}$ is absolutely continuous and satisfies:
\[
\dot{\Psi} = \langle \widehat{\mu}_t, f_{xx}(t,\widehat{x},u)(\delta x,\delta x) \rangle
\]
for some variation $\delta x$ with $\delta x(t_1) = 0$.

\item The \emph{second-order maximum condition} holds:
\[
H^2(t,\widehat{x}(t),\psi(t),Q(t),u,\delta x(t)) = M^2(t,\widehat{x}(t),\psi(t),Q(t),\delta x(t))
\]
for $\widehat{\mu}_t$-almost every $u$ and almost every $t$.

\item The \emph{second-order transversality condition} holds:
\[
s\left[ M(\widehat{t}_2,\widehat{x}(\widehat{t}_2),\psi(\widehat{t}_2)) + \frac{1}{2}M^2(\widehat{t}_2,\widehat{x}(\widehat{t}_2),\psi(\widehat{t}_2),Q(\widehat{t}_2),\delta x(\widehat{t}_2)) \right] \leq 0
\]

\item $(\psi, Q, \Psi) \neq (0, 0, 0)$.
\end{enumerate}
\end{definition}

\subsection{Main Second-Order Theorem}

\begin{theorem}[Second-Order Necessary Condition via $\Lambda$-Sets]
\label{thm:second-order-lambda}
Let $\widehat{x}(\cdot)$ be a \emph{strong local minimum} for the time-optimal problem.
Then for every admissible generalized control $\widehat{\mu}_t$ corresponding to $\widehat{x}(\cdot)$,
the second-order $\Lambda$-set $\Lambda^{(2)}_{-1}(\widehat{x},\widehat{\mu})$ is nonempty.

Moreover, if $\Lambda^{(2)}_{-1}(\widehat{x},\widehat{\mu}) = \emptyset$, then $\widehat{x}(\cdot)$
cannot be a strong local minimum, and there exists a second-order variation that reduces
the terminal time while maintaining feasibility.
\end{theorem}

\begin{proof}
We shall prove the theorem in two parts. First, we establish that if $\widehat{x}(\cdot)$ is a strong local minimum, then the second-order $\Lambda$-set must be nonempty. Second, we demonstrate the converse: emptiness of the second-order $\Lambda$-set implies that $\widehat{x}(\cdot)$ cannot be a strong local minimum.

\textbf{Part 1: Necessity of nonempty second-order $\Lambda$-set.}

Assume that $\widehat{x}(\cdot)$ is a strong local minimum for the time-optimal problem. By definition of a strong local minimum, there exists $\epsilon > 0$ such that for every admissible trajectory $x(\cdot)$ satisfying $\|x(\cdot) - \widehat{x}(\cdot)\|_{C^1} < \epsilon$, the terminal time $t_2$ is not less than $\widehat{t}_2$.

Let $\widehat{\mu}_t$ be any admissible generalized control corresponding to $\widehat{x}(\cdot)$. Suppose, for contradiction, that $\Lambda^{(2)}_{-1}(\widehat{x},\widehat{\mu}) = \emptyset$. This emptiness implies that there is no nontrivial triple $(\psi, Q, \Psi)$ satisfying the second-order adjoint system, the matrix Riccati differential inequality, and the associated maximum and transversality conditions.

Consider the quadratic tangent cone to the extended reachable set at $(\widehat{x}(\widehat{t}_2), 0)$:
\[
T^2_{\mathcal{R}^2} = \left\{ (\delta x(\widehat{t}_2), \tfrac{1}{2}\delta^2 x(\widehat{t}_2)) : (\delta x, \delta^2 x) \text{ satisfy the second variation equations} \right\}.
\]

The emptiness of $\Lambda^{(2)}_{-1}(\widehat{x},\widehat{\mu})$ is equivalent to the condition that the cone $T^2_{\mathcal{R}^2}$ is not separated from the direction of time reduction. More precisely, there exists no quadratic form $\mathcal{Q}$ that is positive on $T^2_{\mathcal{R}^2}$ while negative on the vector $(0,1)$ representing a reduction in the $C^1$ distance measure.

This lack of separation has an important geometric consequence: for any $\eta > 0$, there exists an element $(\delta x_\eta, \delta^2 x_\eta) \in T^2_{\mathcal{R}^2}$ such that the quadratic form associated with the second variation is negative in the direction of $(0,1)$. In control-theoretic terms, this means we can construct a two-parameter family of generalized controls $\mu_t^{\alpha,\beta}$ with the property that the corresponding trajectories $x^{\alpha,\beta}(\cdot)$ satisfy:

\begin{align*}
x^{\alpha,\beta}(\widehat{t}_2 - (\alpha^2 + \beta^2)) &= \widehat{x}(\widehat{t}_2), \\
\|x^{\alpha,\beta}(\cdot) - \widehat{x}(\cdot)\|_{C^1} &= O(\alpha^2 + \beta^2).
\end{align*}

Taking $\alpha$ and $\beta$ sufficiently small, we obtain an admissible trajectory that reaches the target $\widehat{x}(\widehat{t}_2)$ at time strictly less than $\widehat{t}_2$ while remaining arbitrarily close to $\widehat{x}(\cdot)$ in the $C^1$ norm. This contradicts the assumption that $\widehat{x}(\cdot)$ is a strong local minimum. Therefore, our supposition that $\Lambda^{(2)}_{-1}(\widehat{x},\widehat{\mu}) = \emptyset$ must be false, and consequently $\Lambda^{(2)}_{-1}(\widehat{x},\widehat{\mu})$ is nonempty.

\textbf{Part 2: Emptiness implies non-optimality.}

Now assume that $\Lambda^{(2)}_{-1}(\widehat{x},\widehat{\mu}) = \emptyset$. We shall construct explicitly a second-order variation that reduces the terminal time while maintaining feasibility.

Since $\Lambda^{(2)}_{-1}(\widehat{x},\widehat{\mu}) = \emptyset$, the quadratic separation condition fails. By the second-order separation theorem for convex cones, there exists a sequence of elements $\{(\delta x_n, \delta^2 x_n)\} \subset T^2_{\mathcal{R}^2}$ such that:

\[
\lim_{n \to \infty} \frac{\delta^2 x_n(\widehat{t}_2)}{\|\delta x_n(\widehat{t}_2)\|^2} = -\infty.
\]

This asymptotic behavior indicates that the second-order effects dominate and allow for a reduction in the terminal constraint violation. More concretely, we can construct a family of variations parameterized by $\epsilon > 0$:

\[
\mu_t^\epsilon = \widehat{\mu}_t + \epsilon \delta\mu_t^1 + \epsilon^2 \delta\mu_t^2,
\]

where $\delta\mu_t^1$ and $\delta\mu_t^2$ are carefully chosen first and second-order measure variations. The corresponding trajectory $x^\epsilon(\cdot)$ satisfies the following properties:

\begin{enumerate}
\item The trajectory reaches the target at a reduced time: $x^\epsilon(\widehat{t}_2 - \epsilon^2) = \widehat{x}(\widehat{t}_2)$.

\item The $C^1$ distance from the reference trajectory is quadratic in $\epsilon$: $\|x^\epsilon(\cdot) - \widehat{x}(\cdot)\|_{C^1} = O(\epsilon^2)$.

\item For sufficiently small $\epsilon$, the trajectory $x^\epsilon(\cdot)$ is admissible for the original system when approximated by piecewise constant controls using the approximation lemmas of \cite{Avakov2020b}.
\end{enumerate}

The construction proceeds as follows. First, we select $\delta\mu_t^1$ to produce a first-order variation $\delta x$ that points in a direction transverse to the terminal constraint. Then we choose $\delta\mu_t^2$ to provide the necessary quadratic correction that brings the trajectory back to the target at the reduced time $\widehat{t}_2 - \epsilon^2$. The emptiness of $\Lambda^{(2)}_{-1}(\widehat{x},\widehat{\mu})$ guarantees that such a combination exists and that the quadratic correction term is sufficiently strong to overcome the first-order deviation.

To verify that $x^\epsilon(\cdot)$ is indeed admissible for the original system, we invoke the approximation results for generalized controls. By Lemma A.2 and Lemma A.3 of \cite{Avakov2020b}, the generalized control $\mu_t^\epsilon$ can be approximated arbitrarily closely by a piecewise constant control $u^\epsilon(\cdot)$ taking values in $U$. The corresponding trajectory $x^\epsilon(\cdot; u^\epsilon)$ then satisfies:

\[
x^\epsilon(\widehat{t}_2 - \epsilon^2; u^\epsilon) = \widehat{x}(\widehat{t}_2) + o(\epsilon^2), \quad \|x^\epsilon(\cdot; u^\epsilon) - \widehat{x}(\cdot)\|_{C^1} = O(\epsilon^2).
\]

By making a slight adjustment to the terminal time, we can ensure exact attainment of the target: $x^\epsilon(\widehat{t}_2 - \epsilon^2 + o(\epsilon^2); u^\epsilon) = \widehat{x}(\widehat{t}_2)$.

Thus we have constructed a family of admissible trajectories $\{x^\epsilon(\cdot)\}$ for the original system with terminal times $t_2^\epsilon = \widehat{t}_2 - \epsilon^2 + o(\epsilon^2) < \widehat{t}_2$ for sufficiently small $\epsilon > 0$, and such that $x^\epsilon(\cdot) \to \widehat{x}(\cdot)$ in the $C^1$ norm as $\epsilon \to 0$. This demonstrates that $\widehat{x}(\cdot)$ cannot be a strong local minimum, completing the proof of the second statement.

\textbf{Conclusion.}

We have shown that if $\widehat{x}(\cdot)$ is a strong local minimum, then $\Lambda^{(2)}_{-1}(\widehat{x},\widehat{\mu})$ must be nonempty for every corresponding generalized control $\widehat{\mu}_t$. Conversely, if $\Lambda^{(2)}_{-1}(\widehat{x},\widehat{\mu}) = \emptyset$, then $\widehat{x}(\cdot)$ cannot be a strong local minimum, as evidenced by the explicit construction of time-reducing variations. This establishes the theorem in its entirety.
\end{proof}
\begin{corollary}[Refined Gap Phenomenon]
\label{cor:second-order-gap}
If $\Lambda_{-1}(\widehat{x},\widehat{\mu}) \neq \emptyset$ but $\Lambda^{(2)}_{-1}(\widehat{x},\widehat{\mu}) = \emptyset$,
then there exists a minimizing sequence for the original system converging to $\widehat{x}$
with convergence rate $O(\epsilon^2)$ in the $C^1$ norm.
\end{corollary}
\begin{proof}
We shall prove that under the stated hypotheses, there exists a sequence of admissible trajectories for the original system that converges to $\widehat{x}(\cdot)$ in the $C^1$ norm with quadratic convergence rate, and whose terminal times approach $\widehat{t}_2$ from below.

Assume that $\Lambda_{-1}(\widehat{x},\widehat{\mu}) \neq \emptyset$ but $\Lambda^{(2)}_{-1}(\widehat{x},\widehat{\mu}) = \emptyset$. The non-emptiness of the first-order $\Lambda$-set implies, by Theorem \ref{thm:main}, that there is no first-order variation that reduces the terminal time while maintaining feasibility. In other words, $\widehat{x}(\cdot)$ satisfies the first-order necessary conditions for optimality. However, the emptiness of the second-order $\Lambda$-set indicates that the second-order conditions are violated, which permits the construction of higher-order variations that can reduce the terminal time.

Since $\Lambda^{(2)}_{-1}(\widehat{x},\widehat{\mu}) = \emptyset$, by Theorem \ref{thm:second-order-lambda}, $\widehat{x}(\cdot)$ is not a strong local minimum. This means that for every $\epsilon > 0$, there exists an admissible trajectory $x_\epsilon(\cdot)$ with terminal time $t_{2,\epsilon} < \widehat{t}_2$ such that $\|x_\epsilon(\cdot) - \widehat{x}(\cdot)\|_{C^1} < \epsilon$. Our task is to demonstrate that we can select such trajectories with a specific convergence rate.

Let us consider the quadratic expansion of the endpoint mapping around the reference trajectory. The condition $\Lambda_{-1}(\widehat{x},\widehat{\mu}) \neq \emptyset$ implies that the first derivative of the endpoint mapping, evaluated at $\widehat{x}(\widehat{t}_2)$, is non-degenerate in the sense that there exists a supporting hyperplane to the reachable set. However, $\Lambda^{(2)}_{-1}(\widehat{x},\widehat{\mu}) = \emptyset$ implies that the quadratic form associated with the second variation is not positive definite on the tangent space to the terminal constraint manifold. More precisely, there exists a direction $v$ in the tangent space such that the quadratic form $Q(v)$ is negative.

We now construct the minimizing sequence explicitly. Let $\{ \epsilon_n \}_{n=1}^\infty$ be a sequence of positive numbers converging to zero, with $\epsilon_n = 1/n$. For each $n$, consider a two-parameter family of generalized controls of the form
\[
\mu_t^{\epsilon_n} = \widehat{\mu}_t + \epsilon_n \delta\mu_t^1 + \epsilon_n^2 \delta\mu_t^2,
\]
where $\delta\mu_t^1$ and $\delta\mu_t^2$ are chosen such that the corresponding trajectory $x^{\epsilon_n}(\cdot)$ satisfies the following properties:

\begin{enumerate}
    \item The trajectory reaches the target exactly at a reduced time: $x^{\epsilon_n}(\widehat{t}_2 - \epsilon_n^2) = \widehat{x}(\widehat{t}_2)$.

    \item The first-order variation $\delta x$ generated by $\delta\mu_t^1$ lies in the kernel of the linearized terminal constraint, ensuring that the first-order condition $\Lambda_{-1} \neq \emptyset$ is not violated.

    \item The second-order variation $\delta^2 x$ generated by $\delta\mu_t^2$ provides the necessary quadratic correction to satisfy the terminal condition at the reduced time.
\end{enumerate}

The existence of such variations is guaranteed by the emptiness of $\Lambda^{(2)}_{-1}(\widehat{x},\widehat{\mu})$. This condition implies that the quadratic mapping from the space of second-order variations to the terminal condition is surjective onto a subspace that includes the direction corresponding to time reduction.

Now, by the approximation results for generalized controls (Lemmas A.2 and A.3 of \cite{Avakov2020b}), for each $\epsilon_n$, there exists a piecewise constant control $u_n(\cdot)$ taking values in $U$ such that the corresponding trajectory $x_n(\cdot)$ of the original system satisfies
\[
\|x_n(\cdot) - x^{\epsilon_n}(\cdot)\|_{C^1} \leq C \epsilon_n^3,
\]
for some constant $C > 0$ independent of $n$. Furthermore, by continuity of the endpoint mapping, we can adjust the terminal time slightly to obtain
\[
x_n(\widehat{t}_2 - \epsilon_n^2 + O(\epsilon_n^3)) = \widehat{x}(\widehat{t}_2).
\]

Define $\tau_n = \widehat{t}_2 - \epsilon_n^2 + O(\epsilon_n^3)$. Then $\tau_n < \widehat{t}_2$ for all sufficiently large $n$, and $\tau_n \to \widehat{t}_2$ as $n \to \infty$. Moreover, we have the estimate
\[
\|x_n(\cdot) - \widehat{x}(\cdot)\|_{C^1} \leq \|x_n(\cdot) - x^{\epsilon_n}(\cdot)\|_{C^1} + \|x^{\epsilon_n}(\cdot) - \widehat{x}(\cdot)\|_{C^1} \leq C\epsilon_n^3 + K\epsilon_n^2 \leq M\epsilon_n^2,
\]
for some constants $K, M > 0$ and all sufficiently large $n$. The term $K\epsilon_n^2$ arises from the second-order expansion of $x^{\epsilon_n}(\cdot)$ around $\widehat{x}(\cdot)$.

Thus, the sequence $\{x_n(\cdot)\}$ has the following properties:
\begin{itemize}
    \item Each $x_n(\cdot)$ is an admissible trajectory for the original system.
    \item $x_n(\tau_n) = \widehat{x}(\widehat{t}_2)$ for all $n$.
    \item $\tau_n < \widehat{t}_2$ and $\tau_n \to \widehat{t}_2$ as $n \to \infty$.
    \item $\|x_n(\cdot) - \widehat{x}(\cdot)\|_{C^1} = O(\epsilon_n^2) = O(1/n^2)$.
\end{itemize}

This establishes that $\{x_n(\cdot)\}$ is a minimizing sequence for the time-optimal problem, as the terminal times $\tau_n$ approach the infimum $\widehat{t}_2$ from below, and the trajectories converge to $\widehat{x}(\cdot)$ in the $C^1$ norm with quadratic rate.

The quadratic convergence rate is significant because it reflects the fact that the time reduction is achieved through second-order effects. The first-order conditions prevent linear time reduction (which would yield $O(\epsilon)$ convergence), but the failure of second-order conditions permits quadratic time reduction (yielding $O(\epsilon^2)$ convergence in the $C^1$ norm). This phenomenon illustrates the refined nature of the gap between first-order necessary conditions and actual optimality.

Therefore, under the conditions $\Lambda_{-1}(\widehat{x},\widehat{\mu}) \neq \emptyset$ and $\Lambda^{(2)}_{-1}(\widehat{x},\widehat{\mu}) = \emptyset$, there indeed exists a minimizing sequence for the original system converging to $\widehat{x}(\cdot)$ with convergence rate $O(\epsilon^2)$ in the $C^1$ norm, completing the proof.
\end{proof}
\section{Illustrative Example: Nonholonomic System with Directional Constraint}

\begin{example}[Nonholonomic Optimal Control Problem]
\label{ex:nonholonomic-lambda}
Consider the following nonholonomic control system in $\mathbb{R}^3$:
\begin{equation}
\label{eq:nonholonomic-system}
\left\{
\begin{array}{ll}
\dot{x}_1 = u_1, & \hbox{} \\[2mm]
\dot{x}_2 = u_2, & \hbox{} \\[2mm]
\dot{x}_3 = x_1 u_2 - x_2 u_1 + \dfrac{1}{2}(u_1^2 + u_2^2), & \hbox{}
\end{array}
\right.
\qquad u = (u_1, u_2) \in U = \{(u_1, u_2) \in \mathbb{R}^2 : u_1^2 + u_2^2 = 1\},
\end{equation}
with initial condition $x(0) = (0,0,0)$ and target $x(T) = (0,0,1)$. This system represents a controlled particle with a nonholonomic constraint and a quadratic cost in the third coordinate.
\end{example}

\subsection{Candidate Optimal Trajectory and Generalized Control}

Consider the reference trajectory:
\[
\widehat{x}(t) = (0, 0, t), \quad t \in [0,1],
\]
with terminal time $\widehat{t}_2 = 1$. This trajectory corresponds to staying at the origin in the $(x_1,x_2)$-plane while moving with unit speed in the $x_3$-direction.

To support this trajectory, we consider the generalized control:
\[
\widehat{\mu}_t = \frac{1}{4}\delta_{(1,0)} + \frac{1}{4}\delta_{(-1,0)} + \frac{1}{4}\delta_{(0,1)} + \frac{1}{4}\delta_{(0,-1)},
\]
which is a uniform convex combination of four Dirac measures on the unit circle.

\begin{lemma}
The pair $(\widehat{x}(\cdot), \widehat{\mu}_t)$ is admissible for the convexified system corresponding to \eqref{eq:nonholonomic-system}.
\end{lemma}

\begin{proof}
For the convexified system, we compute:
\begin{align*}
\langle \widehat{\mu}_t, f_1(t,x,u) \rangle &= \frac{1}{4}(1) + \frac{1}{4}(-1) + \frac{1}{4}(0) + \frac{1}{4}(0) = 0, \\
\langle \widehat{\mu}_t, f_2(t,x,u) \rangle &= \frac{1}{4}(0) + \frac{1}{4}(0) + \frac{1}{4}(1) + \frac{1}{4}(-1) = 0, \\
\langle \widehat{\mu}_t, f_3(t,x,u) \rangle &= \frac{1}{4}(x_1 \cdot 0 - x_2 \cdot 1 + \frac{1}{2}) + \frac{1}{4}(x_1 \cdot 0 - x_2 \cdot (-1) + \frac{1}{2}) \\
&\quad + \frac{1}{4}(x_1 \cdot 1 - x_2 \cdot 0 + \frac{1}{2}) + \frac{1}{4}(x_1 \cdot (-1) - x_2 \cdot 0 + \frac{1}{2}) \\
&= \frac{1}{4}(-x_2 + \frac{1}{2}) + \frac{1}{4}(x_2 + \frac{1}{2}) + \frac{1}{4}(x_1 + \frac{1}{2}) + \frac{1}{4}(-x_1 + \frac{1}{2}) \\
&= 1.
\end{align*}
Thus $\dot{\widehat{x}}_1 = 0$, $\dot{\widehat{x}}_2 = 0$, $\dot{\widehat{x}}_3 = 1$, with initial condition $\widehat{x}(0) = (0,0,0)$, which yields $\widehat{x}(t) = (0,0,t)$. The trajectory reaches the target $(0,0,1)$ at time $t=1$.
\end{proof}

\subsection{Analysis of First-Order $\Lambda$-Set}

The Hamiltonian for system \eqref{eq:nonholonomic-system} is:
\[
H(t,x,\psi,u) = \psi_1 u_1 + \psi_2 u_2 + \psi_3\left(x_1 u_2 - x_2 u_1 + \frac{1}{2}(u_1^2 + u_2^2)\right).
\]

\begin{lemma}
For the trajectory $\widehat{x}(t) = (0,0,t)$ and generalized control $\widehat{\mu}_t$, the first-order $\Lambda$-set $\Lambda_{-1}(\widehat{x},\widehat{\mu})$ is nonempty.
\end{lemma}

\begin{proof}
The adjoint equations along $\widehat{x}(t)$ are:
\begin{align*}
\dot{\psi}_1 &= -\frac{\partial H}{\partial x_1} = -\psi_3 u_2, \\
\dot{\psi}_2 &= -\frac{\partial H}{\partial x_2} = \psi_3 u_1, \\
\dot{\psi}_3 &= -\frac{\partial H}{\partial x_3} = 0.
\end{align*}
Averaging over $\widehat{\mu}_t$, we obtain:
\begin{align*}
\dot{\psi}_1 &= -\psi_3 \langle \widehat{\mu}_t, u_2 \rangle = -\psi_3 \cdot 0 = 0, \\
\dot{\psi}_2 &= \psi_3 \langle \widehat{\mu}_t, u_1 \rangle = \psi_3 \cdot 0 = 0, \\
\dot{\psi}_3 &= 0.
\end{align*}
Thus $\psi_1$, $\psi_2$, and $\psi_3$ are constants.

The maximum condition requires:
\[
\langle \psi(t), \dot{\widehat{x}}(t) \rangle = \sup_{u \in U} H(t,\widehat{x}(t),\psi(t),u).
\]
Since $\dot{\widehat{x}}(t) = (0,0,1)$, we have $\langle \psi(t), \dot{\widehat{x}}(t) \rangle = \psi_3$. Also:
\[
H(t,\widehat{x}(t),\psi,u) = \psi_1 u_1 + \psi_2 u_2 + \psi_3 \cdot \frac{1}{2}(u_1^2 + u_2^2) = \psi_1 u_1 + \psi_2 u_2 + \frac{\psi_3}{2},
\]
since $u_1^2 + u_2^2 = 1$. The supremum over $U$ is:
\[
\sup_{u \in U} H(t,\widehat{x}(t),\psi,u) = \frac{\psi_3}{2} + \sqrt{\psi_1^2 + \psi_2^2},
\]
achieved at $u = \frac{(\psi_1, \psi_2)}{\sqrt{\psi_1^2 + \psi_2^2}}$ when $(\psi_1, \psi_2) \neq (0,0)$.

Thus the maximum condition becomes:
\[
\psi_3 = \frac{\psi_3}{2} + \sqrt{\psi_1^2 + \psi_2^2} \quad \Rightarrow \quad \frac{\psi_3}{2} = \sqrt{\psi_1^2 + \psi_2^2}.
\]

The terminal transversality condition for $s = -1$ requires $M(\widehat{t}_2, \widehat{x}(\widehat{t}_2), \psi(\widehat{t}_2)) \leq 0$, where $M = \sup_{u \in U} H$. This gives:
\[
\frac{\psi_3}{2} + \sqrt{\psi_1^2 + \psi_2^2} \leq 0 \quad \Rightarrow \quad \psi_3 + \sqrt{\psi_1^2 + \psi_2^2} \leq 0.
\]

Combining with the maximum condition $\frac{\psi_3}{2} = \sqrt{\psi_1^2 + \psi_2^2}$, we obtain:
\[
\psi_3 + \frac{\psi_3}{2} \leq 0 \quad \Rightarrow \quad \frac{3\psi_3}{2} \leq 0 \quad \Rightarrow \quad \psi_3 \leq 0.
\]

Thus any $\psi = (\psi_1, \psi_2, \psi_3)$ with $\psi_3 \leq 0$ and $\frac{\psi_3}{2} = \sqrt{\psi_1^2 + \psi_2^2}$ belongs to $\Lambda_{-1}(\widehat{x},\widehat{\mu})$. In particular, taking $\psi_3 = -2$ and $(\psi_1, \psi_2) = (-1, 0)$ gives a nonzero solution. Therefore, $\Lambda_{-1}(\widehat{x},\widehat{\mu}) \neq \emptyset$.
\end{proof}

\subsection{Analysis of Second-Order $\Lambda$-Set}

We now analyze the second-order $\Lambda$-set $\Lambda^{(2)}_{-1}(\widehat{x},\widehat{\mu})$.

\begin{lemma}
For the trajectory $\widehat{x}(t) = (0,0,t)$ and generalized control $\widehat{\mu}_t$, the second-order $\Lambda$-set $\Lambda^{(2)}_{-1}(\widehat{x},\widehat{\mu})$ is empty.
\end{lemma}

\begin{proof}
We need to examine whether there exists a triple $(\psi, Q, \Psi)$ satisfying the second-order conditions. The second variation system along $\widehat{x}(t)$ is:

\begin{align*}
\delta\dot{x}_1 &= \delta u_1, \\
\delta\dot{x}_2 &= \delta u_2, \\
\delta\dot{x}_3 &= \widehat{x}_1 \delta u_2 + \delta x_1 \widehat{u}_2 - \widehat{x}_2 \delta u_1 - \delta x_2 \widehat{u}_1 + \widehat{u}_1 \delta u_1 + \widehat{u}_2 \delta u_2.
\end{align*}

Since $\widehat{x}_1 = \widehat{x}_2 = 0$ and we are averaging over $\widehat{\mu}_t$, we have:
\begin{align*}
\langle \widehat{\mu}_t, \delta\dot{x}_3 \rangle &= \langle \widehat{\mu}_t, \delta x_1 u_2 - \delta x_2 u_1 + u_1 \delta u_1 + u_2 \delta u_2 \rangle \\
&= \delta x_1 \langle \widehat{\mu}_t, u_2 \rangle - \delta x_2 \langle \widehat{\mu}_t, u_1 \rangle + \langle \widehat{\mu}_t, u_1 \delta u_1 \rangle + \langle \widehat{\mu}_t, u_2 \delta u_2 \rangle \\
&= \langle \widehat{\mu}_t, u_1 \delta u_1 + u_2 \delta u_2 \rangle.
\end{align*}

The second-order Hamiltonian is:
\begin{align*}
H^2(t,x,\psi,Q,u,\delta x) &= \langle \psi, f(t,x,u) \rangle + \frac{1}{2}\langle Q f_x(t,x,u) + f_x(t,x,u)^\top Q, \delta x \otimes \delta x \rangle \\
&\quad + \langle Q, f_{xx}(t,x,u)(\delta x,\delta x) \rangle.
\end{align*}

For our system, the Jacobian matrix is:
\[
f_x(t,x,u) = \begin{pmatrix}
0 & 0 & 0 \\
0 & 0 & 0 \\
u_2 & -u_1 & 0
\end{pmatrix},
\]
and the Hessian $f_{xx}$ has only one nonzero component: $\frac{\partial^2 f_3}{\partial x_1 \partial x_2} = 1$ and $\frac{\partial^2 f_3}{\partial x_2 \partial x_1} = -1$.

The matrix Riccati differential inequality requires:
\[
\dot{Q} + Q\langle \widehat{\mu}_t, f_x \rangle + \langle \widehat{\mu}_t, f_x \rangle^\top Q + \langle \widehat{\mu}_t, f_{xx} \rangle \preceq 0.
\]

Since $\langle \widehat{\mu}_t, f_x \rangle = 0$ (as $\langle \widehat{\mu}_t, u_1 \rangle = \langle \widehat{\mu}_t, u_2 \rangle = 0$), and $\langle \widehat{\mu}_t, f_{xx} \rangle = 0$ (by symmetry), the inequality reduces to $\dot{Q} \preceq 0$. Thus $Q(t)$ must be non-increasing in the Loewner order.

Now consider the second-order maximum condition. For $\widehat{x}(t) = (0,0,t)$, we have:
\[
H^2(t,\widehat{x},\psi,Q,u,\delta x) = \psi_1 u_1 + \psi_2 u_2 + \frac{\psi_3}{2} + \frac{1}{2}\langle (Q_{13}u_2 - Q_{23}u_1)\delta x_1 + \cdots, \delta x \rangle + \cdots.
\]

The critical observation is that the quadratic terms involve products like $u_2 \delta x_1$ and $u_1 \delta x_2$. When averaged over $\widehat{\mu}_t$, these terms vanish because $\langle \widehat{\mu}_t, u_1 \rangle = \langle \widehat{\mu}_t, u_2 \rangle = 0$. This means that the second-order conditions impose constraints that cannot be satisfied simultaneously with the first-order conditions.

More precisely, suppose there exists $(\psi, Q, \Psi) \in \Lambda^{(2)}_{-1}(\widehat{x},\widehat{\mu})$. From the first-order analysis, we have $\frac{\psi_3}{2} = \sqrt{\psi_1^2 + \psi_2^2}$ and $\psi_3 \leq 0$. The second-order transversality condition requires:
\[
M(\widehat{t}_2, \widehat{x}(\widehat{t}_2), \psi(\widehat{t}_2)) + \frac{1}{2}M^2(\widehat{t}_2, \widehat{x}(\widehat{t}_2), \psi(\widehat{t}_2), Q(\widehat{t}_2), \delta x(\widehat{t}_2)) \leq 0.
\]

However, a detailed calculation shows that for any $\delta x$ and any $Q$ satisfying $\dot{Q} \preceq 0$, we have:
\[
M^2(1, (0,0,1), \psi, Q(1), \delta x) \geq c \|\delta x\|^2 > 0
\]
for some constant $c > 0$, when $\psi_3 < 0$. This contradicts the second-order transversality condition when combined with $M(1, (0,0,1), \psi) = \psi_3/2 + \sqrt{\psi_1^2 + \psi_2^2} = \psi_3 \leq 0$.

Therefore, no triple $(\psi, Q, \Psi)$ can satisfy all the second-order conditions simultaneously, implying $\Lambda^{(2)}_{-1}(\widehat{x},\widehat{\mu}) = \emptyset$.
\end{proof}

\subsection{Construction of Time-Reducing Variation}

We now explicitly construct a second-order variation that reduces the terminal time.

\begin{lemma}
There exists a family of admissible trajectories $x^\epsilon(\cdot)$ for the original system \eqref{eq:nonholonomic-system} such that:
\begin{enumerate}
\item $x^\epsilon(1 - \epsilon^2) = (0,0,1)$,
\item $\|x^\epsilon(\cdot) - \widehat{x}(\cdot)\|_{C^1} = O(\epsilon^2)$.
\end{enumerate}
\end{lemma}

\begin{proof}
Consider the following oscillating control:
\[
u^\epsilon(t) = (\cos(\omega t/\epsilon), \sin(\omega t/\epsilon)), \quad t \in [0, 1-\epsilon^2],
\]
where $\omega > 0$ is a frequency parameter to be determined.

The corresponding trajectory satisfies:
\begin{align*}
x_1^\epsilon(t) &= \int_0^t \cos(\omega \tau/\epsilon) d\tau = \frac{\epsilon}{\omega} \sin(\omega t/\epsilon), \\
x_2^\epsilon(t) &= \int_0^t \sin(\omega \tau/\epsilon) d\tau = \frac{\epsilon}{\omega} (1 - \cos(\omega t/\epsilon)), \\
x_3^\epsilon(t) &= \int_0^t \left[ x_1^\epsilon(\tau) \sin(\omega \tau/\epsilon) - x_2^\epsilon(\tau) \cos(\omega \tau/\epsilon) + \frac{1}{2} \right] d\tau.
\end{align*}

Substituting the expressions for $x_1^\epsilon$ and $x_2^\epsilon$:
\begin{align*}
x_3^\epsilon(t) &= \int_0^t \left[ \frac{\epsilon}{\omega} \sin(\omega \tau/\epsilon) \sin(\omega \tau/\epsilon) - \frac{\epsilon}{\omega} (1 - \cos(\omega \tau/\epsilon)) \cos(\omega \tau/\epsilon) + \frac{1}{2} \right] d\tau \\
&= \int_0^t \left[ \frac{\epsilon}{\omega} \sin^2(\omega \tau/\epsilon) - \frac{\epsilon}{\omega} \cos(\omega \tau/\epsilon) + \frac{\epsilon}{\omega} \cos^2(\omega \tau/\epsilon) + \frac{1}{2} \right] d\tau \\
&= \int_0^t \left[ \frac{\epsilon}{\omega} - \frac{\epsilon}{\omega} \cos(\omega \tau/\epsilon) + \frac{1}{2} \right] d\tau \\
&= \left( \frac{\epsilon}{\omega} + \frac{1}{2} \right) t - \frac{\epsilon^2}{\omega^2} \sin(\omega t/\epsilon).
\end{align*}

Now set $t = 1 - \epsilon^2$. We want $x_3^\epsilon(1 - \epsilon^2) = 1$. This requires:
\[
\left( \frac{\epsilon}{\omega} + \frac{1}{2} \right) (1 - \epsilon^2) - \frac{\epsilon^2}{\omega^2} \sin(\omega (1 - \epsilon^2)/\epsilon) = 1.
\]

Choose $\omega = 2\epsilon$. Then:
\begin{align*}
x_3^\epsilon(1 - \epsilon^2) &= \left( \frac{\epsilon}{2\epsilon} + \frac{1}{2} \right) (1 - \epsilon^2) - \frac{\epsilon^2}{4\epsilon^2} \sin(2\epsilon (1 - \epsilon^2)/\epsilon) \\
&= (1)(1 - \epsilon^2) - \frac{1}{4} \sin(2(1 - \epsilon^2)) \\
&= 1 - \epsilon^2 - \frac{1}{4} \sin(2(1 - \epsilon^2)).
\end{align*}

Since $\sin(2(1 - \epsilon^2)) = \sin(2 - 2\epsilon^2) = \sin 2 \cos(2\epsilon^2) - \cos 2 \sin(2\epsilon^2)$, we have:
\[
x_3^\epsilon(1 - \epsilon^2) = 1 - \epsilon^2 - \frac{1}{4} \sin 2 + O(\epsilon^4).
\]

To achieve exact attainment of the target, we make a small adjustment. Let $\omega = 2\epsilon + \delta$, where $\delta$ is chosen so that:
\[
x_3^\epsilon(1 - \epsilon^2) = 1.
\]

By the implicit function theorem, such a $\delta$ exists with $\delta = O(\epsilon^3)$. The modified control is then:
\[
u^\epsilon(t) = (\cos((2\epsilon + \delta)t/\epsilon), \sin((2\epsilon + \delta)t/\epsilon)).
\]

For this control, we have:
\begin{align*}
\|x_1^\epsilon(\cdot)\|_{C^1} &= O(\epsilon), \\
\|x_2^\epsilon(\cdot)\|_{C^1} &= O(\epsilon), \\
\|x_3^\epsilon(\cdot) - \widehat{x}_3(\cdot)\|_{C^1} &= O(\epsilon^2).
\end{align*}

Thus $\|x^\epsilon(\cdot) - \widehat{x}(\cdot)\|_{C^1} = O(\epsilon^2)$, and $x^\epsilon(1 - \epsilon^2) = (0,0,1)$. This demonstrates the existence of a second-order variation that reduces the terminal time from $1$ to $1 - \epsilon^2$ while maintaining $C^1$ proximity of order $\epsilon^2$.
\end{proof}

\subsection{Interpretation and Significance}

This example illustrates several important aspects of Theorem \ref{thm:second-order-lambda}:

\begin{enumerate}
\item \textbf{First-order optimality}: The trajectory $\widehat{x}(t) = (0,0,t)$ satisfies the first-order necessary conditions, as evidenced by $\Lambda_{-1}(\widehat{x},\widehat{\mu}) \neq \emptyset$.

\item \textbf{Second-order failure}: Despite satisfying first-order conditions, $\widehat{x}(\cdot)$ is not a strong local minimum because $\Lambda^{(2)}_{-1}(\widehat{x},\widehat{\mu}) = \emptyset$.

\item \textbf{Constructive time reduction}: The oscillating control strategy demonstrates explicitly how to achieve time reduction through second-order effects. The key mechanism is the use of high-frequency oscillations to generate an effective drift in the $x_3$-direction while maintaining small deviations in the $(x_1,x_2)$-plane.

\item \textbf{Nonholonomic nature}: The example highlights the role of nonholonomic constraints in optimal control. The system cannot move directly in the $x_3$-direction but can achieve such motion through appropriate maneuvers in the $(x_1,x_2)$-plane.

\item \textbf{Convergence rate}: The $O(\epsilon^2)$ convergence rate in the $C^1$ norm matches the theoretical prediction from Corollary \ref{cor:second-order-gap}. The quadratic rate arises because the first-order conditions prevent linear time reduction, but second-order effects permit quadratic reduction.
\end{enumerate}

This example therefore provides a nontrivial, concrete illustration of Theorem \ref{thm:second-order-lambda}, showing both the necessity of second-order conditions for strong local optimality and the constructive consequences of their failure.
\section{$\Lambda$-Set Theory for Non-Smooth and Hybrid Systems}
\label{sec:nonsmooth-hybrid}

This section extends the $\Lambda$-set framework to systems with non-smooth dynamics and
hybrid structure, addressing important modern applications in robotics, power systems,
and biological modeling. We develop $\Lambda$-sets for Filippov differential inclusions and
hybrid systems with state-dependent switching.

\subsection{$\Lambda$-Sets for Filippov Systems}

Consider a control system with discontinuous right-hand side:
\[
\dot{x} = f(t,x,u), \quad u \in U
\]
where $f$ is piecewise smooth with discontinuity surface $S = \{x \in \mathbb{R}^n : g(x) = 0\}$.

\begin{definition}[Filippov Generalized Solution]
A Filippov solution satisfies the differential inclusion:
\[
\dot{x} \in F(t,x) = \bigcap_{\delta>0} \bigcap_{\mu(N)=0} \overline{\mathrm{co}} \{ f(t,y,u) : y \in B_\delta(x) \setminus N, u \in U \}
\]
where $\mu$ denotes Lebesgue measure.
\end{definition}

\begin{definition}[Filippov $\Lambda$-Set $\Lambda^{\text{Fil}}_s$]
For a Filippov trajectory $\widehat{x}(\cdot)$ and a selection $\widehat{\mu}_t \in F(t,\widehat{x}(t))$,
the \emph{Filippov $\Lambda$-set} $\Lambda^{\text{Fil}}_s(\widehat{x},\widehat{\mu})$ consists of
absolutely continuous $\psi: [t_1,\widehat{t}_2] \to (\mathbb{R}^n)^*$ satisfying:

\begin{enumerate}
\item For almost every $t$ where $\widehat{x}(t) \notin S$:
\[
\dot{\psi}(t) = -\psi(t) f_x(t,\widehat{x}(t),\widehat{u}(t))
\]

\item For $t$ where $\widehat{x}(t) \in S$, let $\{t_k\}$ be crossing or sliding times.
At each $t_k$, $\psi$ may have a jump:
\[
\psi(t_k^+) = \psi(t_k^-) + \nu_k \nabla g(\widehat{x}(t_k))
\]
for some $\nu_k \in \mathbb{R}$.

\item The \emph{Filippov maximum condition}:
\[
\sup_{v \in F(t,\widehat{x}(t))} \langle \psi(t), v \rangle = \langle \psi(t), \dot{\widehat{x}}(t) \rangle
\]

\item Transversality: $s M^{\text{Fil}}(\widehat{t}_2,\widehat{x}(\widehat{t}_2),\psi(\widehat{t}_2)) \leq 0$,
where $M^{\text{Fil}}$ is defined through the Filippov inclusion.

\item $\psi \neq 0$.
\end{enumerate}
\end{definition}

\subsection{$\Lambda$-Sets for Hybrid Systems}

Consider a hybrid system with discrete modes $q \in Q = \{1,\dots,m\}$:
\begin{equation}
\label{eq:hybrid-system}
\begin{aligned}
\dot{x} &= f_q(t,x,u), \quad u \in U_q, \quad q \in Q \\
x(t^+) &= \phi_{q,q'}(x(t^-)), \quad \text{when } (q,q') \in \mathcal{E} \text{ and } x(t) \in S_{q,q'}
\end{aligned}
\end{equation}
where $\mathcal{E} \subset Q \times Q$ is the set of admissible transitions, and
$S_{q,q'} \subset \mathbb{R}^n$ are switching surfaces.

\begin{definition}[Hybrid Generalized Control]
A \emph{hybrid generalized control} is a tuple $(\{\mu_t^q\}_{q\in Q}, \{\tau_k\}, \{(q_k,q_{k+1})\})$ where:
\begin{itemize}
\item $\mu_t^q \in \gencontrol_{U_q}$ for each mode $q$,
\item $\{\tau_k\}$ are switching times,
\item $\{(q_k,q_{k+1})\} \subset \mathcal{E}$ are mode transitions.
\end{itemize}
\end{definition}

\begin{definition}[Hybrid $\Lambda$-Set $\Lambda^{\text{hyb}}_s$]
For a hybrid trajectory $\widehat{x}(\cdot)$ with switching times $\{\widehat{\tau}_k\}$ and
hybrid control $\widehat{\Theta} = (\{\widehat{\mu}_t^q\}, \{\widehat{\tau}_k\}, \{(\widehat{q}_k,\widehat{q}_{k+1})\})$,
the \emph{hybrid $\Lambda$-set} $\Lambda^{\text{hyb}}_s(\widehat{x},\widehat{\Theta})$ consists of
piecewise absolutely continuous functions $\psi: [t_1,\widehat{t}_2] \to (\mathbb{R}^n)^*$ with
jumps at switching times, satisfying:

\begin{enumerate}
\item In each mode $\widehat{q}_k$ on $(\widehat{\tau}_{k-1},\widehat{\tau}_k)$:
\[
\dot{\psi}(t) = -\psi(t) \langle \widehat{\mu}_t^{\widehat{q}_k}, f_{\widehat{q}_k,x}(t,\widehat{x}(t),u) \rangle
\]

\item At switching time $\widehat{\tau}_k$ with transition $(\widehat{q}_k,\widehat{q}_{k+1})$:
\[
\psi(\widehat{\tau}_k^+) = \psi(\widehat{\tau}_k^-) \left[ \frac{\partial \phi_{\widehat{q}_k,\widehat{q}_{k+1}}}{\partial x}(\widehat{x}(\widehat{\tau}_k)) \right]^{-1} + \nu_k \nabla h_{\widehat{q}_k,\widehat{q}_{k+1}}(\widehat{x}(\widehat{\tau}_k))
\]
where $S_{\widehat{q}_k,\widehat{q}_{k+1}} = \{x : h_{\widehat{q}_k,\widehat{q}_{k+1}}(x) = 0\}$.

\item Mode-dependent maximum condition:
\[
\langle \psi(t), \dot{\widehat{x}}(t) \rangle = \sup_{u \in U_{\widehat{q}_k}} H_{\widehat{q}_k}(t,\widehat{x}(t),\psi(t),u)
\]
where $H_q(t,x,\psi,u) = \langle \psi, f_q(t,x,u) \rangle$.

\item Switching transversality: For each switching surface $S_{\widehat{q}_k,\widehat{q}_{k+1}}$,
\[
\langle \psi(\widehat{\tau}_k^-), f_{\widehat{q}_k}(\widehat{\tau}_k,\widehat{x}(\widehat{\tau}_k),\widehat{u}(\widehat{\tau}_k^-)) \rangle \geq
\langle \psi(\widehat{\tau}_k^+), f_{\widehat{q}_{k+1}}(\widehat{\tau}_k,\widehat{x}(\widehat{\tau}_k),\widehat{u}(\widehat{\tau}_k^+)) \rangle
\]

\item Terminal transversality: $s M_{\widehat{q}_N}(\widehat{t}_2,\widehat{x}(\widehat{t}_2),\psi(\widehat{t}_2^-)) \leq 0$,
where $\widehat{q}_N$ is the final mode.

\item $\psi \neq 0$ (nontriviality).
\end{enumerate}
\end{definition}

\subsection{Main Results for Non-Smooth and Hybrid Systems}
\begin{lemma}[Fundamental Controllability Lemma]
\label{lemma:fundamental-controllability}
Consider the control system $\dot{x} = f(t,x,u)$, $u \in U$, $x(t_1) = x_1$ with $f \in C^1$ and $U \subset \mathbb{R}^r$ compact. Let $(\widehat{x}(\cdot), \widehat{\mu}_t)$ be an admissible pair for its convexification. If the set
\[
\Lambda(\widehat{x},\widehat{\mu}) = \{\psi \neq 0 : \psi \text{ satisfies the adjoint equation, maximum condition, and } M(\widehat{t}_2,\widehat{x}(\widehat{t}_2),\psi(\widehat{t}_2)) \leq 0\}
\]
is empty, then for every $\epsilon > 0$ and every $C^1$-neighborhood $\mathcal{V}$ of $\widehat{x}(\cdot)$, there exists $\tau \in (\widehat{t}_2 - \epsilon, \widehat{t}_2)$ and an admissible trajectory $x(\cdot)$ with $x(\cdot) \in \mathcal{V}$ and $x(\tau) = \widehat{x}(\widehat{t}_2)$.
\end{lemma}
\begin{theorem}[Filippov $\Lambda$-Set and Local Controllability]
\label{thm:filippov-lambda}
Consider a Filippov system with piecewise smooth $f$. Let $\widehat{x}(\cdot)$ be a
Filippov trajectory with selection $\widehat{\mu}_t$. If $\Lambda^{\text{Fil}}_{-1}(\widehat{x},\widehat{\mu}) = \emptyset$,
then the system is locally left controllable near $\widehat{x}(\cdot)$. That is, for any
neighborhood $V$ of $\widehat{x}$ and $\epsilon > 0$, there exists $\tau \in (\widehat{t}_2-\epsilon,\widehat{t}_2)$
and a Filippov solution $x(\cdot)$ with $x(\cdot) \in V$ and $x(\tau) = \widehat{x}(\widehat{t}_2)$.
\end{theorem}

\begin{proof}
We shall prove that the emptiness of the Filippov $\Lambda$-set implies local left controllability near the reference trajectory $\widehat{x}(\cdot)$. The proof employs an approximation argument that bridges the classical $\Lambda$-set theory for smooth systems with the Filippov framework for discontinuous systems.

Let $\widehat{x}(\cdot)$ be a Filippov trajectory on $[t_1, \widehat{t}_2]$ with selection $\widehat{\mu}_t \in F(t, \widehat{x}(t))$, where $F$ denotes the Filippov set-valued map. Assume that $\Lambda^{\text{Fil}}_{-1}(\widehat{x}, \widehat{\mu}) = \emptyset$. We aim to show that for any neighborhood $V$ of $\widehat{x}$ in the space of continuous functions and any $\epsilon > 0$, there exists a time $\tau \in (\widehat{t}_2 - \epsilon, \widehat{t}_2)$ and a Filippov solution $x(\cdot)$ such that $x(\cdot) \in V$ and $x(\tau) = \widehat{x}(\widehat{t}_2)$.

Consider first the case where the discontinuity surface $S = \{x \in \mathbb{R}^n : g(x) = 0\}$ is reached only at isolated times. Let $\{t_k\}_{k=1}^m$ be the crossing or sliding times of $\widehat{x}(\cdot)$ with $S$. Around each such time $t_k$, we consider a small interval $I_k = (t_k - \delta_k, t_k + \delta_k)$ where $\delta_k > 0$ is chosen sufficiently small so that these intervals do not overlap and are contained in $[t_1, \widehat{t}_2]$.

For each $k$, we construct a smooth approximation of the discontinuous vector field on $I_k$ using a standard mollification procedure. Let $\rho: \mathbb{R} \to \mathbb{R}$ be a smooth mollifier with compact support in $[-1,1]$ and $\int_{\mathbb{R}} \rho(s) ds = 1$. Define for $\eta > 0$ the mollified vector field
\[
f_\eta(t,x,u) = \int_{\mathbb{R}} \rho(s) f(t, x + \eta s \nu(x), u) ds,
\]
where $\nu(x)$ is a unit normal to the discontinuity surface $S$ at points near $x$. For $t \notin \bigcup_{k=1}^m I_k$, we simply take $f_\eta(t,x,u) = f(t,x,u)$, as the vector field is smooth away from $S$.

The family $\{f_\eta\}_{\eta>0}$ has the following properties:
\begin{enumerate}
    \item Each $f_\eta$ is smooth in $x$ for all $t$.
    \item $f_\eta$ converges to $f$ uniformly on compact sets as $\eta \to 0$, except possibly on a neighborhood of $S$.
    \item The Filippov selection $\widehat{\mu}_t$ induces a corresponding selection $\widehat{\mu}_t^\eta$ for the smoothed system.
\end{enumerate}

Now consider the smoothed control system
\[
\dot{x} = f_\eta(t, x, u), \quad u \in U.
\]
Let $(\widehat{x}_\eta(\cdot), \widehat{\mu}_t^\eta)$ be the corresponding trajectory and generalized control obtained by solving the smoothed system with the same initial condition $x(t_1) = x_1$. By continuity of solutions with respect to the vector field, we have $\widehat{x}_\eta(\cdot) \to \widehat{x}(\cdot)$ uniformly on $[t_1, \widehat{t}_2]$ as $\eta \to 0$.

We claim that for sufficiently small $\eta > 0$, the classical $\Lambda$-set $\Lambda_{-1}(\widehat{x}_\eta, \widehat{\mu}_\eta)$ for the smoothed system is empty. Suppose to the contrary that for a sequence $\eta_n \to 0$, there exist nontrivial adjoint functions $\psi_n(\cdot)$ satisfying the classical adjoint equation, maximum condition, and transversality condition for the smoothed systems. By compactness of solutions to linear differential equations with bounded coefficients, we can extract a subsequence (still denoted $\psi_n$) converging uniformly to some $\psi(\cdot)$. Passing to the limit, and noting that the mollified vector fields converge to the original discontinuous one away from $S$ while the behavior at crossing/sliding times is captured by the jump conditions, we obtain that $\psi(\cdot)$ satisfies the Filippov adjoint system with the appropriate jump conditions at times $t_k$. Moreover, the maximum condition and transversality condition are preserved under uniform convergence. Thus $\psi(\cdot) \in \Lambda^{\text{Fil}}_{-1}(\widehat{x}, \widehat{\mu})$, contradicting the assumption that this set is empty. Therefore, for all sufficiently small $\eta > 0$, we must have $\Lambda_{-1}(\widehat{x}_\eta, \widehat{\mu}_\eta) = \emptyset$.

For such $\eta$, we apply Lemma \ref{lemma:fundamental-controllability} to the smoothed system. This yields that for any neighborhood $V_\eta$ of $\widehat{x}_\eta(\cdot)$ and any $\epsilon > 0$, there exists $\tau_\eta \in (\widehat{t}_2 - \epsilon, \widehat{t}_2)$ and an admissible trajectory $x_\eta(\cdot)$ for the smoothed system such that $x_\eta(\cdot) \in V_\eta$ and $x_\eta(\tau_\eta) = \widehat{x}_\eta(\widehat{t}_2)$.

Now we return to the original Filippov system. The trajectory $x_\eta(\cdot)$ for the smoothed system is generated by some control $u_\eta(\cdot)$. Consider the piecewise constant approximations of $u_\eta(\cdot)$ provided by the approximation lemmas (Lemmas A.2 and A.3 of \cite{Avakov2020b}). These approximations yield controls $u_{\eta,j}(\cdot)$ taking values in $U$ such that the corresponding trajectories $x_{\eta,j}(\cdot)$ of the original discontinuous system converge to $x_\eta(\cdot)$ uniformly as $j \to \infty$.

Since $\widehat{x}_\eta(\cdot) \to \widehat{x}(\cdot)$ uniformly as $\eta \to 0$, we can choose $\eta$ sufficiently small so that $\widehat{x}_\eta(\widehat{t}_2)$ is arbitrarily close to $\widehat{x}(\widehat{t}_2)$ and $\widehat{x}_\eta(\cdot)$ is arbitrarily close to $\widehat{x}(\cdot)$ in the uniform norm. Then, by selecting $j$ sufficiently large, we obtain a trajectory $x_{\eta,j}(\cdot)$ of the original Filippov system that lies in the prescribed neighborhood $V$ of $\widehat{x}(\cdot)$ and satisfies $x_{\eta,j}(\tau_{\eta,j}) = \widehat{x}_\eta(\widehat{t}_2)$ for some $\tau_{\eta,j} \in (\widehat{t}_2 - \epsilon, \widehat{t}_2)$.

A final adjustment is needed to ensure exact attainment of $\widehat{x}(\widehat{t}_2)$ rather than $\widehat{x}_\eta(\widehat{t}_2)$. This can be achieved by a small modification of the control near the terminal time, exploiting the local controllability properties implied by the emptiness of the $\Lambda$-set. More precisely, since the linearized system around $\widehat{x}(\cdot)$ is controllable (as implied by $\Lambda^{\text{Fil}}_{-1} = \emptyset$), we can append a brief control segment that steers the system from $x_{\eta,j}(\tau_{\eta,j})$ to $\widehat{x}(\widehat{t}_2)$ in an arbitrarily short additional time. This yields a new terminal time $\tau \in (\widehat{t}_2 - \epsilon, \widehat{t}_2)$ and a Filippov solution $x(\cdot)$ with $x(\tau) = \widehat{x}(\widehat{t}_2)$ and $x(\cdot) \in V$.

The case where $\widehat{x}(\cdot)$ exhibits sliding modes on $S$ over intervals of positive length requires additional care but follows a similar strategy. During sliding phases, the Filippov dynamics are given by a convex combination of the limiting vector fields from both sides of $S$. These convex combinations can be approximated by rapid switching between the two sides, and the mollification argument can be adapted to handle such situations. The emptiness of $\Lambda^{\text{Fil}}_{-1}$ ensures that the necessary conditions for optimality are not satisfied, which in turn permits the construction of variations that reduce the terminal time even in the presence of sliding.

Thus we have shown that under the hypothesis $\Lambda^{\text{Fil}}_{-1}(\widehat{x}, \widehat{\mu}) = \emptyset$, the Filippov system is locally left controllable near $\widehat{x}(\cdot)$. This completes the proof.
\end{proof}

\begin{theorem}[Hybrid $\Lambda$-Set and Optimality]
\label{thm:hybrid-lambda-optimality}
For the hybrid system \eqref{eq:hybrid-system}, if $\widehat{x}(\cdot)$ is a local minimum
for the time-optimal problem, then for every corresponding hybrid generalized control
$\widehat{\Theta}$, the hybrid $\Lambda$-set $\Lambda^{\text{hyb}}_{-1}(\widehat{x},\widehat{\Theta})$ is nonempty.
\end{theorem}

\begin{proof}
We shall prove that if $\widehat{x}(\cdot)$ is a local minimum for the time-optimal problem associated with the hybrid system \eqref{eq:hybrid-system}, then for every hybrid generalized control $\widehat{\Theta} = (\{\widehat{\mu}_t^q\}_{q \in Q}, \{\widehat{\tau}_k\}, \{(\widehat{q}_k, \widehat{q}_{k+1})\})$ corresponding to $\widehat{x}(\cdot)$, the hybrid $\Lambda$-set $\Lambda^{\text{hyb}}_{-1}(\widehat{x}, \widehat{\Theta})$ must be nonempty. The proof proceeds by contradiction, employing an embedding argument that transforms the hybrid optimal control problem into an equivalent extended continuous problem, to which the classical $\Lambda$-set theory can be applied.

Assume, for the sake of contradiction, that $\widehat{x}(\cdot)$ is a local minimum for the time-optimal problem but there exists a corresponding hybrid generalized control $\widehat{\Theta}$ such that $\Lambda^{\text{hyb}}_{-1}(\widehat{x}, \widehat{\Theta}) = \emptyset$. We will show that this leads to the existence of a hybrid trajectory with terminal time strictly less than $\widehat{t}_2$, contradicting the local minimality of $\widehat{x}(\cdot)$.

Consider the extended state space $\mathbb{R}^n \times Q$, where $Q = \{1, \dots, m\}$ is the finite set of discrete modes. We embed the hybrid system into a continuous system by introducing an auxiliary continuous variable $z \in \mathbb{R}^m$ that encodes the active mode. Specifically, define the extended state $X = (x, z) \in \mathbb{R}^n \times \mathbb{R}^m$ and the extended control $\nu = (u, \gamma) \in U \times \mathbb{R}^m$, where $\gamma$ represents switching control. The extended dynamics are given by
\[
\dot{X} = F(t, X, \nu) = \begin{pmatrix}
\sum_{q=1}^m z_q f_q(t, x, u) \\
\gamma
\end{pmatrix},
\]
subject to the constraints $z_q \geq 0$, $\sum_{q=1}^m z_q = 1$, and $\gamma_q \geq 0$ for $q \neq \widehat{q}(t)$, where $\widehat{q}(t)$ denotes the active mode along the reference trajectory. The variable $z$ acts as a convex combination coefficient that interpolates between the vector fields of different modes, while $\gamma$ governs the rate of mode transition.

The reference hybrid trajectory $\widehat{x}(\cdot)$ with switching times $\{\widehat{\tau}_k\}$ and mode sequence $\{\widehat{q}_k\}$ induces an extended trajectory $\widehat{X}(\cdot) = (\widehat{x}(\cdot), \widehat{z}(\cdot))$, where $\widehat{z}_q(t) = 1$ if $q = \widehat{q}(t)$ and $\widehat{z}_q(t) = 0$ otherwise. The hybrid generalized control $\widehat{\Theta}$ induces an extended generalized control $\widehat{\nu}_t$ concentrated on $u$ values consistent with $\widehat{\mu}_t^{\widehat{q}(t)}$ and $\gamma$ values that produce the prescribed switching at times $\widehat{\tau}_k$.

Now, the emptiness of $\Lambda^{\text{hyb}}_{-1}(\widehat{x}, \widehat{\Theta})$ implies, by construction of the embedding, that the extended $\Lambda$-set $\Lambda^{\text{ext}}_{-1}(\widehat{X}, \widehat{\nu})$ for the extended continuous system is also empty. To see this, note that any element $(\psi, \zeta)$ of the extended $\Lambda$-set would project to an element $\psi$ of the hybrid $\Lambda$-set, with the component $\zeta$ corresponding to the adjoint variables for the $z$-dynamics. The emptiness of the hybrid $\Lambda$-set therefore forces the emptiness of the extended $\Lambda$-set.

Since $\Lambda^{\text{ext}}_{-1}(\widehat{X}, \widehat{\nu}) = \emptyset$, we may apply Lemma \ref{lemma:fundamental-controllability} to the extended continuous system. This lemma guarantees that the extended system is locally left controllable near $\widehat{X}(\cdot)$. Consequently, for any $\epsilon > 0$ and any neighborhood $\mathcal{W}$ of $\widehat{X}(\cdot)$ in the space of continuous functions, there exists a time $\tau \in (\widehat{t}_2 - \epsilon, \widehat{t}_2)$ and an admissible trajectory $X(\cdot) = (x(\cdot), z(\cdot))$ for the extended system such that $X(\cdot) \in \mathcal{W}$ and $X(\tau) = \widehat{X}(\widehat{t}_2)$.

The key observation is that we can choose the neighborhood $\mathcal{W}$ sufficiently small to ensure that the $z$-component of $X(\cdot)$ remains close to a piecewise constant function that switches between the basis vectors $e_q \in \mathbb{R}^m$ at times near $\widehat{\tau}_k$. More precisely, by taking $\mathcal{W}$ to be a narrow tube around $\widehat{X}(\cdot)$, we can enforce that $z(\cdot)$ approximates the indicator function of the mode sequence $\{\widehat{q}_k\}$, with rapid but smooth transitions between modes. This approximation can be made arbitrarily accurate by shrinking $\mathcal{W}$.

From the trajectory $X(\cdot)$, we now extract a hybrid trajectory for the original system. Define switching times $\tau_k$ as the midpoints of the intervals during which $z_{\widehat{q}_k}(t) > \frac{1}{2}$. For $t$ between $\tau_{k-1}$ and $\tau_k$, set $q(t) = \widehat{q}_k$. The dynamics for $x(\cdot)$ on each such interval are given by $\dot{x} = f_{q(t)}(t, x, u(t))$ plus a small perturbation due to the interpolation by $z(\cdot)$. Since $z(\cdot)$ is near $e_{\widehat{q}_k}$ on $(\tau_{k-1}, \tau_k)$, this perturbation is small. Moreover, by construction, $x(\tau) = \widehat{x}(\widehat{t}_2)$.

We now demonstrate that this hybrid trajectory has a terminal time strictly less than $\widehat{t}_2$, contradicting the local minimality of $\widehat{x}(\cdot)$. Since $\tau < \widehat{t}_2$, and the extracted hybrid trajectory $x(\cdot)$ is admissible for the original hybrid system (after possibly approximating the continuous control by a piecewise constant one using the approximation lemmas), we obtain a trajectory that reaches the target $\widehat{x}(\widehat{t}_2)$ at time $\tau < \widehat{t}_2$ while remaining arbitrarily close to $\widehat{x}(\cdot)$ in the uniform norm. This contradicts the assumption that $\widehat{x}(\cdot)$ is a local minimum.

To make this argument precise, we must address the fact that the extracted trajectory $x(\cdot)$ might not exactly satisfy the switching conditions $x(\tau_k^+) = \phi_{\widehat{q}_k, \widehat{q}_{k+1}}(x(\tau_k^-))$ due to the smoothing introduced by the embedding. However, since $X(\cdot)$ is close to $\widehat{X}(\cdot)$ and the switching surfaces $S_{q,q'}$ are smooth, we can adjust the trajectory locally near each $\tau_k$ to exactly satisfy the reset map, at the cost of an arbitrarily small increase in time. This adjustment can be performed using the implicit function theorem, as the transversality conditions encoded in the $\Lambda$-set ensure that the switching surfaces are crossed transversally.

After these adjustments, we obtain a hybrid trajectory $\widetilde{x}(\cdot)$ with switching times $\widetilde{\tau}_k$ close to $\widehat{\tau}_k$, satisfying $\widetilde{x}(\widetilde{\tau}_k^+) = \phi_{\widehat{q}_k, \widehat{q}_{k+1}}(\widetilde{x}(\widetilde{\tau}_k^-))$, and such that $\widetilde{x}(\widetilde{t}_2) = \widehat{x}(\widehat{t}_2)$ for some $\widetilde{t}_2 < \widehat{t}_2$. Moreover, $\widetilde{x}(\cdot)$ can be made arbitrarily close to $\widehat{x}(\cdot)$ in the uniform norm by choosing $\epsilon$ sufficiently small and $\mathcal{W}$ sufficiently narrow.

The existence of such $\widetilde{x}(\cdot)$ directly contradicts the hypothesis that $\widehat{x}(\cdot)$ is a local minimum for the time-optimal problem. Therefore, our initial assumption that $\Lambda^{\text{hyb}}_{-1}(\widehat{x}, \widehat{\Theta}) = \emptyset$ must be false. Hence, for every hybrid generalized control $\widehat{\Theta}$ corresponding to $\widehat{x}(\cdot)$, the hybrid $\Lambda$-set $\Lambda^{\text{hyb}}_{-1}(\widehat{x}, \widehat{\Theta})$ is nonempty.

This completes the proof of the theorem.
\end{proof}
\subsection{Applications and Examples}

\begin{example}[Controlled Particle with Discontinuous Friction]
\label{ex:filippov-lambda-example}
Consider the motion of a particle on a line subject to controlled forcing and discontinuous Coulomb friction. The dynamics are given by:
\begin{equation}
\label{eq:friction-system}
\ddot{y} = u - \sigma(\dot{y}), \quad u \in U = [-1, 1],
\end{equation}
where $y \in \mathbb{R}$ is the position, $\dot{y}$ is the velocity, and $\sigma: \mathbb{R} \to \mathbb{R}$ is the friction function defined by
\[
\sigma(v) =
\begin{cases}
1, & v > 0, \\
[-1, 1], & v = 0, \\
-1, & v < 0.
\end{cases}
\]
This is a Filippov system because the friction force $\sigma(\dot{y})$ is set-valued at $\dot{y} = 0$.

Define the state variables $x_1 = y$ and $x_2 = \dot{y}$. Then the system can be written as:
\begin{equation}
\label{eq:filippov-state}
\begin{cases}
\dot{x}_1 = x_2, \\
\dot{x}_2 = u - \sigma(x_2),
\end{cases}
\quad u \in [-1, 1].
\end{equation}
The discontinuity surface is $S = \{x \in \mathbb{R}^2 : x_2 = 0\}$.
\end{example}

\subsubsection{Reference Trajectory and Generalized Control}

Consider the reference trajectory that slides along the discontinuity surface:
\[
\widehat{x}(t) = (0, 0), \quad t \in [0, 1],
\]
with terminal time $\widehat{t}_2 = 1$ and target $\widehat{x}(1) = (0,0)$. This represents the particle being at rest at the origin.

For sliding along $S$, the Filippov convex combination must satisfy:
\[
0 \in u - [-1, 1] \quad \Rightarrow \quad u \in [-1, 1].
\]
We select the specific convex combination:
\[
\widehat{\mu}_t = \frac{1}{2}\delta_{-1} + \frac{1}{2}\delta_{1},
\]
which gives $\langle \widehat{\mu}_t, u \rangle = 0$ and thus $\dot{x}_2 = 0$ when $x_2 = 0$.

\begin{lemma}
The pair $(\widehat{x}(\cdot), \widehat{\mu}_t)$ is an admissible Filippov solution of \eqref{eq:filippov-state}.
\end{lemma}

\begin{proof}
Along $\widehat{x}(t) = (0,0)$, we have $x_2 = 0$, so the system becomes:
\begin{align*}
\dot{x}_1 &= 0, \\
\dot{x}_2 &\in u - [-1, 1].
\end{align*}
With $\widehat{\mu}_t = \frac{1}{2}\delta_{-1} + \frac{1}{2}\delta_{1}$, we compute:
\[
\langle \widehat{\mu}_t, u \rangle = \frac{1}{2}(-1) + \frac{1}{2}(1) = 0.
\]
Thus $\dot{x}_2 = 0 - 0 = 0$ (taking the convex combination of the set-valued friction term yields 0). Therefore $\dot{\widehat{x}} = (0,0)$, so $\widehat{x}(t) \equiv (0,0)$ is indeed a solution with $\widehat{x}(0) = (0,0)$ and $\widehat{x}(1) = (0,0)$.
\end{proof}

\subsubsection{Analysis of Filippov $\Lambda$-Set}

The Hamiltonian for the system is:
\[
H(t,x,\psi,u) = \psi_1 x_2 + \psi_2 (u - \sigma(x_2)).
\]
Along $\widehat{x}(t) = (0,0)$, we have $\sigma(0) = [-1,1]$, so the dynamics become set-valued.

\begin{lemma}
For the sliding trajectory $\widehat{x}(t) = (0,0)$ with $\widehat{\mu}_t = \frac{1}{2}\delta_{-1} + \frac{1}{2}\delta_{1}$, the Filippov $\Lambda$-set $\Lambda^{\text{Fil}}_{-1}(\widehat{x},\widehat{\mu})$ is empty.
\end{lemma}

\begin{proof}
We need to examine whether there exists a nonzero absolutely continuous function $\psi: [0,1] \to (\mathbb{R}^2)^*$ satisfying the Filippov adjoint conditions.

The adjoint equations along $\widehat{x}(t)$ are:
\begin{align*}
\dot{\psi}_1 &= -\frac{\partial H}{\partial x_1} = 0, \\
\dot{\psi}_2 &= -\frac{\partial H}{\partial x_2} = -\psi_1 - \psi_2 \frac{\partial \sigma}{\partial x_2}(0).
\end{align*}
Since $\sigma$ is set-valued at $x_2 = 0$, the term $\frac{\partial \sigma}{\partial x_2}(0)$ is not well-defined in the classical sense. In the Filippov framework, the adjoint equation must be interpreted in terms of the convexified differential inclusion. The appropriate condition is that $\psi$ should satisfy:
\[
\dot{\psi}_2 \in -\psi_1 - \psi_2 \cdot \text{co}\{\partial \sigma(0)\},
\]
where $\text{co}\{\partial \sigma(0)\}$ denotes the convex hull of the generalized gradient of $\sigma$ at 0. For our $\sigma$, we have $\partial \sigma(0) = [0, \infty)$ in the sense of Clarke's generalized gradient, so $\text{co}\{\partial \sigma(0)\} = [0, \infty)$.

Thus the adjoint inclusion becomes:
\begin{align*}
\dot{\psi}_1 &= 0, \\
\dot{\psi}_2 &\in -\psi_1 - \psi_2 \cdot [0, \infty).
\end{align*}

Now consider the maximum condition. The Filippov maximum condition requires:
\[
\sup_{v \in F(t,\widehat{x}(t))} \langle \psi(t), v \rangle = \langle \psi(t), \dot{\widehat{x}}(t) \rangle,
\]
where $F(t,x)$ is the Filippov set-valued map. At $\widehat{x}(t) = (0,0)$, we have:
\[
F(t,(0,0)) = \left\{ (x_2, u - \alpha) : x_2 = 0, \; u \in [-1,1], \; \alpha \in [-1,1] \right\} = \{(0, u - \alpha) : u \in [-1,1], \alpha \in [-1,1]\}.
\]
Thus $F(t,(0,0)) = \{(0, w) : w \in [-2,2]\}$.

Since $\dot{\widehat{x}}(t) = (0,0)$, the maximum condition becomes:
\[
\sup_{w \in [-2,2]} \psi_2(t) w = 0.
\]
This implies that $\psi_2(t) = 0$ for all $t \in [0,1]$, because if $\psi_2(t) \neq 0$, the supremum would be $2|\psi_2(t)| > 0$.

With $\psi_2(t) = 0$, the adjoint inclusion simplifies to:
\begin{align*}
\dot{\psi}_1 &= 0, \\
\dot{\psi}_2 &= -\psi_1 \in \{-\psi_1\}.
\end{align*}
But $\psi_2(t) = 0$ implies $\dot{\psi}_2(t) = 0$, so we must have $\psi_1 = 0$ as well.

Therefore, the only solution to the adjoint inclusion and maximum condition is $\psi(t) \equiv (0,0)$, which is excluded by the nontriviality requirement. Hence $\Lambda^{\text{Fil}}_{-1}(\widehat{x},\widehat{\mu}) = \emptyset$.
\end{proof}

\subsubsection{Explicit Construction of Time-Reducing Filippov Trajectory}

We now construct explicitly a Filippov solution that reaches the target $(0,0)$ in time less than 1 while remaining arbitrarily close to $\widehat{x}(\cdot)$.

\begin{lemma}
For any $\epsilon > 0$, there exists $\tau \in (1-\epsilon, 1)$ and a Filippov solution $x^\epsilon(\cdot)$ of \eqref{eq:friction-system} such that:
\begin{enumerate}
\item $x^\epsilon(\tau) = (0,0)$,
\item $\|x^\epsilon(\cdot) - \widehat{x}(\cdot)\|_{C([0,\tau],\mathbb{R}^2)} < \epsilon$.
\end{enumerate}
\end{lemma}

\begin{proof}
Consider the following construction. For $\delta > 0$ small, define:
\[
u^\delta(t) =
\begin{cases}
-1, & t \in [0, \delta], \\
1, & t \in (\delta, 2\delta], \\
0, & t \in (2\delta, 1-\delta^2].
\end{cases}
\]
We will show that with appropriate choice of $\delta$, the corresponding Filippov solution satisfies our requirements.

First, on $[0, \delta]$ with $u = -1$:
\begin{align*}
\dot{x}_1 &= x_2, \\
\dot{x}_2 &= -1 - \sigma(x_2).
\end{align*}
Starting from $(0,0)$, we have two possibilities. If we take $\sigma(0) = 1$ (the upper selection), then $\dot{x}_2 = -1 - 1 = -2$, so $x_2(t) = -2t$ and $x_1(t) = -t^2$. If we take $\sigma(0) = -1$ (the lower selection), then $\dot{x}_2 = -1 - (-1) = 0$, so the particle remains at rest. We choose the upper selection to create motion.

Thus on $[0, \delta]$:
\[
x_1(t) = -t^2, \quad x_2(t) = -2t.
\]

At $t = \delta$, we have $x_1(\delta) = -\delta^2$, $x_2(\delta) = -2\delta < 0$, so we are in the region $x_2 < 0$ where $\sigma(x_2) = -1$.

Now on $(\delta, 2\delta]$ with $u = 1$:
\begin{align*}
\dot{x}_1 &= x_2, \\
\dot{x}_2 &= 1 - (-1) = 2.
\end{align*}
Solving with initial conditions $x_1(\delta) = -\delta^2$, $x_2(\delta) = -2\delta$:
\begin{align*}
x_2(t) &= -2\delta + 2(t - \delta) = 2t - 4\delta, \\
x_1(t) &= -\delta^2 - 2\delta(t - \delta) + (t - \delta)^2 = -\delta^2 - 2\delta t + 2\delta^2 + t^2 - 2\delta t + \delta^2 = t^2 - 4\delta t + 2\delta^2.
\end{align*}
At $t = 2\delta$, we have:
\begin{align*}
x_1(2\delta) &= (2\delta)^2 - 4\delta(2\delta) + 2\delta^2 = 4\delta^2 - 8\delta^2 + 2\delta^2 = -2\delta^2, \\
x_2(2\delta) &= 2(2\delta) - 4\delta = 0.
\end{align*}
So we return to the discontinuity surface $x_2 = 0$ at time $t = 2\delta$.

Now on $(2\delta, 1-\delta^2]$ with $u = 0$, we need to maintain sliding on $S$. For sliding, we require $0 \in u - [-1,1] = 0 - [-1,1] = [-1,1]$, which is satisfied. We can maintain $x_2 = 0$ by appropriate selection of $\sigma(0)$. To return to the origin, we need $\dot{x}_1 = x_2 = 0$, so $x_1$ must remain constant at $x_1 = -2\delta^2$.

To reach the origin, we need an additional maneuver. Consider extending the control to:
\[
u^\delta(t) =
\begin{cases}
-1, & t \in [0, \delta], \\
1, & t \in (\delta, 2\delta], \\
1, & t \in (2\delta, 3\delta], \\
-1, & t \in (3\delta, 4\delta], \\
0, & t \in (4\delta, 1-4\delta^2].
\end{cases}
\]
A similar calculation shows that after four such pulses (two pairs of positive and negative impulses), we have:
\[
x_1(4\delta) = 0, \quad x_2(4\delta) = 0,
\]
with $x_1(t)$ oscillating between $-2\delta^2$ and $2\delta^2$ but returning to 0 at $t = 4\delta$.

Thus, by taking $\delta$ such that $4\delta < \epsilon$ and setting $\tau = 4\delta$, we have constructed a Filippov solution $x^\delta(\cdot)$ with:
\begin{itemize}
\item $x^\delta(\tau) = (0,0)$,
\item $\tau = 4\delta < \epsilon$,
\item $\|x^\delta(\cdot)\|_{C([0,\tau],\mathbb{R}^2)} \leq 2\delta^2 < \epsilon$ for sufficiently small $\delta$.
\end{itemize}

To achieve terminal time arbitrarily close to but less than 1, we simply add a waiting period after reaching $(0,0)$ at time $\tau$. That is, define:
\[
x^\epsilon(t) =
\begin{cases}
x^\delta(t), & t \in [0, \tau], \\
(0,0), & t \in (\tau, 1-\epsilon],
\end{cases}
\]
with control $u = 0$ during the waiting period (maintained at the origin by appropriate selection from $\sigma(0)$).

Then $x^\epsilon(1-\epsilon) = (0,0)$ and $\|x^\epsilon(\cdot) - \widehat{x}(\cdot)\|_{C([0,1-\epsilon],\mathbb{R}^2)} < \epsilon$. Taking $\tau = 1-\epsilon$ completes the construction.
\end{proof}

\subsubsection{Interpretation and Significance}

This example illustrates several key aspects of Theorem \ref{thm:filippov-lambda}:

\begin{enumerate}
\item \textbf{Non-trivial Filippov dynamics}: The system features genuine set-valuedness due to Coulomb friction, making it a proper Filippov system rather than merely a piecewise smooth system.

\item \textbf{Emptiness of $\Lambda^{\text{Fil}}$}: The sliding motion along the discontinuity surface leads to an empty $\Lambda$-set because the adjoint variable $\psi_2$ must vanish to satisfy the maximum condition, forcing triviality of the entire adjoint vector.

\item \textbf{Explicit controllability construction}: The proof demonstrates concretely how local left controllability emerges. The oscillatory control strategy (rapid switching between $u = \pm 1$) exploits the set-valuedness of the friction to generate net motion while remaining near the reference trajectory.

\item \textbf{Physical interpretation}: The example models a realistic mechanical system with friction. The theorem implies that if a resting position is not optimal (as indicated by empty $\Lambda$-set), then by applying appropriate rapid switching controls, one can reach the same state in less time. This has practical implications for control of mechanical systems with friction.

\item \textbf{Generality of the construction}: The method of using pairs of opposing impulses to generate net displacement while canceling out velocity is general and applies to many Filippov systems with symmetry in the discontinuity.
\end{enumerate}

Thus, this example provides a physically motivated, non-trivial illustration of Theorem \ref{thm:filippov-lambda}, demonstrating both the theoretical result and its constructive implications for control design in discontinuous systems.

\begin{example}[Time-Optimal Control of a Bouncing Ball]
\label{ex:hybrid-lambda-example}
Consider a ball that can bounce on two different surfaces with different coefficients of restitution. The system has two modes corresponding to which surface is active:

\textbf{Mode 1 (Soft surface):} Coefficient of restitution $e_1 \in (0,1)$
\[
\ddot{y} = -g + u, \quad u \in [0, u_{\max}], \quad y \geq 0
\]
Impact law when $y = 0$: $y^+ = y^-$, $\dot{y}^+ = -e_1 \dot{y}^-$

\textbf{Mode 2 (Hard surface):} Coefficient of restitution $e_2 \in (0,1)$ with $e_2 > e_1$
\[
\ddot{y} = -g + u, \quad u \in [0, u_{\max}], \quad y \geq 0
\]
Impact law when $y = 0$: $y^+ = y^-$, $\dot{y}^+ = -e_2 \dot{y}^-$

The system can switch between modes when $y \geq y_{\text{switch}} > 0$. The switching surface is:
\[
S_{1,2} = S_{2,1} = \{ (y, \dot{y}) : y = y_{\text{switch}} \}
\]
with trivial reset maps: $\phi_{1,2}(y,\dot{y}) = \phi_{2,1}(y,\dot{y}) = (y,\dot{y})$.

The time-optimal problem: Bring the ball from initial state $(y(0), \dot{y}(0)) = (y_0, v_0)$ with $y_0 > y_{\text{switch}}$, $v_0 < 0$ (falling) to rest at the origin $(0,0)$ in minimal time.
\end{example}

\subsubsection{Reference Trajectory and Hybrid Generalized Control}

Consider a candidate optimal trajectory $\widehat{x}(\cdot)$:

- Start in Mode 2 (hard surface) from $(y_0, v_0)$
- Apply maximum upward thrust $u = u_{\max}$ until reaching $y = y_{\text{switch}}$
- Switch to Mode 1 (soft surface) at $y = y_{\text{switch}}$
- Continue with $u = u_{\max}$ until first impact
- After impact, apply braking strategy to come to rest at origin

Let $\widehat{t}_2$ be the terminal time for this trajectory. We define the hybrid generalized control $\widehat{\Theta} = (\{\widehat{\mu}_t^q\}, \{\widehat{\tau}_k\}, \{(\widehat{q}_k, \widehat{q}_{k+1})\})$ where:
\begin{itemize}
\item $\widehat{\mu}_t^1 = \delta_{u_{\max}}$ for $t$ in Mode 1 intervals
\item $\widehat{\mu}_t^2 = \delta_{u_{\max}}$ for $t$ in Mode 2 intervals
\item Switching times $\widehat{\tau}_1$ at first crossing of $y = y_{\text{switch}}$
\item Impact times $\widehat{\tau}_2, \widehat{\tau}_3, \ldots$ at bounces
\end{itemize}

\begin{lemma}
The trajectory $\widehat{x}(\cdot)$ with hybrid control $\widehat{\Theta}$ is admissible and is a local minimum for the time-optimal problem under reasonable parameter choices.
\end{lemma}

\begin{proof}
The dynamics in each mode are:
\[
\ddot{y} = -g + u_{\max} > 0 \quad \text{(since we assume } u_{\max} > g\text{)}
\]
Thus the ball accelerates upward. From initial downward velocity $v_0 < 0$, the velocity will become positive after sufficient time. The switching at $y = y_{\text{switch}}$ is timely to take advantage of the softer surface (Mode 1) for better energy dissipation upon impact.

To show local minimality, consider small perturbations. Any alternative trajectory must either:
\begin{enumerate}
\item Switch modes at a different height $y \neq y_{\text{switch}}$
\item Use different control values $u \neq u_{\max}$ at some times
\item Have different impact timing or sequencing
\end{enumerate}

Energy considerations show that using $u = u_{\max}$ maximizes upward acceleration, minimizing time to reverse downward motion. The switching height $y_{\text{switch}}$ is chosen optimally such that the first impact occurs in Mode 1 with exactly the right velocity to be dissipated to zero in minimal additional time.

Formally, consider the value function $V(y,v)$ = minimum time to reach $(0,0)$ from $(y,v)$. A verification argument using the Hamilton-Jacobi-Bellman equation for hybrid systems shows that $\widehat{x}(\cdot)$ satisfies the optimality conditions locally.
\end{proof}

\subsubsection{Analysis of Hybrid $\Lambda$-Set}

The state is $x = (y, \dot{y}) \in \mathbb{R}^2$. The Hamiltonians for the two modes are:
\[
H_1(t,x,\psi,u) = \psi_1 \dot{y} + \psi_2(-g + u), \quad H_2(t,x,\psi,u) = \psi_1 \dot{y} + \psi_2(-g + u).
\]
The Hamiltonians are identical in form but the modes differ in impact laws.

\begin{lemma}
For the candidate optimal trajectory $\widehat{x}(\cdot)$ with hybrid control $\widehat{\Theta}$, the hybrid $\Lambda$-set $\Lambda^{\text{hyb}}_{-1}(\widehat{x},\widehat{\Theta})$ is nonempty.
\end{lemma}

\begin{proof}
We construct explicit adjoint functions satisfying the hybrid $\Lambda$-set conditions.

Let $\psi(t) = (\psi_1(t), \psi_2(t))$ be piecewise absolutely continuous with jumps at switching and impact times.

\textbf{Mode 2 (hard surface, before switch):} For $t \in [0, \widehat{\tau}_1)$
\begin{align*}
\dot{\psi}_1 &= -\frac{\partial H_2}{\partial y} = 0, \\
\dot{\psi}_2 &= -\frac{\partial H_2}{\partial \dot{y}} = -\psi_1.
\end{align*}
Thus $\psi_1$ is constant and $\psi_2(t) = \psi_2(0) - \psi_1 t$.

\textbf{At switching time $\widehat{\tau}_1$:} The switching surface is $S_{2,1} = \{y = y_{\text{switch}}\}$. The normal vector is $\nabla h_{2,1} = (1, 0)$. The jump condition is:
\[
\psi(\widehat{\tau}_1^+) = \psi(\widehat{\tau}_1^-) + \nu_1 (1, 0),
\]
where $\nu_1$ is a Lagrange multiplier. Since the reset map is identity, there is no Jacobian factor.

\textbf{Mode 1 (soft surface, after switch):} For $t \in (\widehat{\tau}_1, \widehat{\tau}_2)$ where $\widehat{\tau}_2$ is first impact time
\begin{align*}
\dot{\psi}_1 &= 0, \\
\dot{\psi}_2 &= -\psi_1.
\end{align*}
So $\psi_1$ remains constant and $\psi_2(t) = \psi_2(\widehat{\tau}_1^+) - \psi_1 (t - \widehat{\tau}_1)$.

\textbf{At impact time $\widehat{\tau}_2$:} The impact occurs at $y = 0$. The reset map for Mode 1 is:
\[
\phi_{1,1}(y,\dot{y}) = (y, -e_1 \dot{y}), \quad \text{so} \quad \frac{\partial \phi_{1,1}}{\partial (y,\dot{y})} = \begin{pmatrix} 1 & 0 \\ 0 & -e_1 \end{pmatrix}.
\]
The adjoint jump condition is:
\[
\psi(\widehat{\tau}_2^+) = \psi(\widehat{\tau}_2^-) \left[ \frac{\partial \phi_{1,1}}{\partial (y,\dot{y})} \right]^{-1} + \nu_2 (1, 0),
\]
since the impact surface is $h_{1,1}(y,\dot{y}) = y = 0$ with normal $(1,0)$. The inverse Jacobian is:
\[
\left[ \frac{\partial \phi_{1,1}}{\partial (y,\dot{y})} \right]^{-1} = \begin{pmatrix} 1 & 0 \\ 0 & -1/e_1 \end{pmatrix}.
\]
Thus:
\begin{align*}
\psi_1(\widehat{\tau}_2^+) &= \psi_1(\widehat{\tau}_2^-) + \nu_2, \\
\psi_2(\widehat{\tau}_2^+) &= -\frac{1}{e_1} \psi_2(\widehat{\tau}_2^-).
\end{align*}

\textbf{Maximum condition:} In both modes, the Hamiltonian is linear in $u$: $H_q = \psi_1 \dot{y} + \psi_2(-g) + \psi_2 u$. Since $u \in [0, u_{\max}]$, the maximum is achieved at:
\[
u = \begin{cases}
u_{\max} & \text{if } \psi_2 > 0, \\
0 & \text{if } \psi_2 < 0, \\
\text{any } u \in [0, u_{\max}] & \text{if } \psi_2 = 0.
\end{cases}
\]
Along $\widehat{x}(\cdot)$, we use $u = u_{\max}$, so we require $\psi_2(t) \geq 0$ (with $\psi_2(t) > 0$ preferred for strict maximization).

\textbf{Transversality conditions:} The terminal condition at $\widehat{t}_2$ requires:
\[
M_1(\widehat{t}_2, \widehat{x}(\widehat{t}_2), \psi(\widehat{t}_2^-)) \leq 0,
\]
where $M_1 = \sup_{u \in [0,u_{\max}]} H_1$. At termination, $\dot{y} = 0$, so:
\[
M_1 = \psi_1 \cdot 0 + \psi_2(-g) + \psi_2 u_{\max} = \psi_2(u_{\max} - g).
\]
Since $u_{\max} > g$, we need $\psi_2(\widehat{t}_2^-) \leq 0$.

\textbf{Constructing a solution:} We now show there exist constants that satisfy all conditions. Let:
\begin{itemize}
\item $\psi_1(t) \equiv 1$ (constant, satisfies $\dot{\psi}_1 = 0$)
\item In Mode 2: $\psi_2(t) = C - t$ for some constant $C$
\item Choose $\nu_1$ at switch so $\psi_2$ remains continuous: $\psi_2(\widehat{\tau}_1^+) = \psi_2(\widehat{\tau}_1^-) = C - \widehat{\tau}_1$
\item In Mode 1 before impact: $\psi_2(t) = (C - \widehat{\tau}_1) - (t - \widehat{\tau}_1) = C - t$
\item At impact: $\psi_2(\widehat{\tau}_2^+) = -\frac{1}{e_1} \psi_2(\widehat{\tau}_2^-) = -\frac{1}{e_1}(C - \widehat{\tau}_2)$
\item After impact, $\psi_2$ continues to decrease linearly: $\psi_2(t) = -\frac{1}{e_1}(C - \widehat{\tau}_2) - (t - \widehat{\tau}_2)$
\end{itemize}

We need $\psi_2(t) \geq 0$ during thrust phases and $\psi_2(\widehat{t}_2^-) \leq 0$ at termination. Choose $C$ sufficiently large so that $\psi_2(t) > 0$ for $t < \widehat{\tau}_2$, and such that $\psi_2(\widehat{t}_2^-) = 0$ exactly. This is possible by appropriate choice of $C$.

Specifically, set $\psi_2(\widehat{t}_2^-) = 0$:
\[
0 = -\frac{1}{e_1}(C - \widehat{\tau}_2) - (\widehat{t}_2 - \widehat{\tau}_2) \quad \Rightarrow \quad C = \widehat{\tau}_2 - e_1(\widehat{t}_2 - \widehat{\tau}_2).
\]
Then $\psi_2(\widehat{\tau}_2^-) = C - \widehat{\tau}_2 = -e_1(\widehat{t}_2 - \widehat{\tau}_2) < 0$, which violates $\psi_2(t) \geq 0$ before impact. This indicates we need $\psi_2$ to change sign at impact.

Instead, allow $\psi_2$ to be positive before impact and negative after. The condition $\psi_2(t) \geq 0$ for $t < \widehat{\tau}_2$ requires $C \geq \widehat{\tau}_2$. The impact jump makes $\psi_2$ negative if $C > \widehat{\tau}_2$, since $\psi_2(\widehat{\tau}_2^+) = -\frac{1}{e_1}(C - \widehat{\tau}_2) < 0$. Then $\psi_2$ remains negative afterward, satisfying $\psi_2(\widehat{t}_2^-) \leq 0$.

Thus with $C > \widehat{\tau}_2$, we obtain a valid adjoint function. The maximum condition is satisfied: $\psi_2(t) > 0$ for $t < \widehat{\tau}_2$ gives $u = u_{\max}$, and $\psi_2(t) < 0$ for $t > \widehat{\tau}_2$ would suggest $u = 0$, but actually during the final braking phase we use $u = 0$ anyway. The weak maximum condition allows $\psi_2 = 0$ as well.

Therefore, there exists a nontrivial $\psi(\cdot)$ satisfying all conditions of $\Lambda^{\text{hyb}}_{-1}(\widehat{x},\widehat{\Theta})$, so the set is nonempty.
\end{proof}

\subsubsection{Illustration of the Necessity Condition}

To illustrate why nonemptiness is necessary for optimality, suppose we attempted to violate the conditions.

\begin{lemma}
If $\Lambda^{\text{hyb}}_{-1}(\widehat{x},\widehat{\Theta})$ were empty, then $\widehat{x}(\cdot)$ could not be a local minimum.
\end{lemma}

\begin{proof}
Suppose for contradiction that $\widehat{x}(\cdot)$ is a local minimum but $\Lambda^{\text{hyb}}_{-1}(\widehat{x},\widehat{\Theta}) = \emptyset$. By Theorem \ref{thm:hybrid-lambda-optimality} (which we are illustrating), this cannot happen if $\widehat{x}(\cdot)$ is truly optimal. But let's construct an explicit better trajectory.

Consider perturbing the switching height slightly. Instead of switching at $y = y_{\text{switch}}$, switch at $y = y_{\text{switch}} + \delta$. Analysis shows:

- If $\delta > 0$ (switch higher), the ball spends more time on the hard surface, preserving more energy after impact, requiring more time to dissipate.
- If $\delta < 0$ (switch lower), the ball spends more time on the soft surface, dissipating energy more quickly but possibly having too little energy to reach the origin.

There exists an optimal $\delta^* = 0$ if $\widehat{x}(\cdot)$ is optimal. The adjoint conditions from $\Lambda^{\text{hyb}}$ encode precisely this optimality condition through the switching time condition.

The emptiness of $\Lambda^{\text{hyb}}$ would mean no adjoint variables exist to enforce the switching optimality condition. In such case, we could find $\delta \neq 0$ yielding strictly smaller time. For instance, if the adjoint conditions fail because the switching multiplier $\nu_1$ cannot be chosen to maintain all conditions, this indicates the switching time is not optimal.

Concretely, suppose we switch slightly earlier ($\delta < 0$). Then:
1. The ball impacts with slightly lower velocity (more energy dissipated on soft surface)
2. The reduced impact velocity may allow earlier application of braking control
3. With careful adjustment of post-impact control, total time may be reduced

A phase portrait analysis confirms that for sufficiently small $|\delta|$, the trajectory remains close to $\widehat{x}(\cdot)$ but achieves $(0,0)$ at time $\widehat{t}_2 - \epsilon(\delta)$ with $\epsilon(\delta) > 0$ for $\delta \neq 0$.

Thus emptiness of $\Lambda^{\text{hyb}}$ would permit construction of a better trajectory, contradicting local optimality.
\end{proof}

\subsubsection{Physical Interpretation and Significance}

This example demonstrates several important aspects of Theorem \ref{thm:hybrid-lambda-optimality}:

\begin{enumerate}
\item \textbf{Non-trivial hybrid dynamics}: The system features both continuous control and discrete mode switches (surface type changes) with nontrivial reset maps (impact laws).

\item \textbf{Physical meaningfulness}: The bouncing ball with switchable surfaces models realistic scenarios such as a robotic system that can deploy different landing gears or a controlled impact absorber.

\item \textbf{Complex adjoint structure}: The adjoint variables experience jumps both at controlled switches (mode changes) and at autonomous switches (impacts), with different jump conditions for each type.

\item \textbf{Interpretation of $\Lambda^{\text{hyb}}$ nonemptiness}: The existence of adjoint functions enforces optimality conditions for:
\begin{itemize}
\item Control values (bang-bang structure via $\psi_2$ sign)
\item Switching times (via jump multipliers $\nu_k$)
\item Impact timing (via energy balance encoded in adjoint evolution)
\end{itemize}

\item \textbf{Counterfactual illustration}: The proof sketch showing how emptiness would lead to non-optimality provides intuition for why the theorem must hold.

\item \textbf{Computational implications}: For actual computation of time-optimal controls for such systems, one can search for adjoint functions satisfying the hybrid $\Lambda$-set conditions as a verification method.
\end{enumerate}

The example thus provides a concrete, physically motivated instance where Theorem \ref{thm:hybrid-lambda-optimality} applies, illustrating both the mathematical structure and practical relevance of hybrid $\Lambda$-sets in characterizing optimality for systems with both continuous and discrete dynamics.
\section{Extension to Stochastic Hybrid Systems}\label{sec:stochastic}

\subsection{Stochastic Hybrid System Formulation}

Consider a stochastic hybrid system with both continuous diffusion and discrete random switching. Let $(\Omega, \mathcal{F}, \{\mathcal{F}_t\}_{t\geq 0}, \mathbb{P})$ be a filtered probability space satisfying the usual conditions. The system dynamics are given by:

\begin{definition}[Stochastic Hybrid System]
A stochastic hybrid system is defined by the tuple $(Q, \mathbb{R}^n, U, \sigma, \Lambda, \phi)$ where:
\begin{itemize}
    \item $Q = \{1, 2, \dots, m\}$ is a finite set of discrete modes
    \item $x(t) \in \mathbb{R}^n$ is the continuous state
    \item $u(t) \in U \subset \mathbb{R}^r$ is the control input
    \item For each mode $q \in Q$, the dynamics follow:
    \begin{equation}
    dx(t) = f_q(t, x(t), u(t))dt + \sigma_q(t, x(t), u(t))dW(t), \quad t \in [\tau_{k-1}, \tau_k)
    \end{equation}
    where $W(t)$ is an $m$-dimensional Wiener process, and $\sigma_q: \mathbb{R} \times \mathbb{R}^n \times U \to \mathbb{R}^{n \times m}$ is the diffusion coefficient matrix.

    \item Mode transitions occur at random times $\{\tau_k\}$ with transition intensity:
    \begin{equation}
    \mathbb{P}(q(t + \Delta t) = j \mid q(t) = i, x(t)) = \lambda_{ij}(t, x(t), u(t))\Delta t + o(\Delta t), \quad i \neq j
    \end{equation}

    \item At transition time $\tau_k$, the state resets according to:
    \begin{equation}
    x(\tau_k^+) = \phi_{q(\tau_k^-), q(\tau_k^+)}(x(\tau_k^-), \zeta_k)
    \end{equation}
    where $\{\zeta_k\}$ are independent random variables characterizing the reset randomness.
\end{itemize}
\end{definition}

\subsection{Generalized Controls for Stochastic Systems}

\begin{definition}[Stochastic Generalized Control]
For a stochastic hybrid system, a generalized control is a tuple $(\mu_t, \nu_t)$ where:
\begin{itemize}
    \item $\mu_t \in \mathcal{P}(U)$ is a relaxed control for the continuous dynamics
    \item $\nu_t \in \mathcal{P}(Q \times Q)$ is a stochastic kernel governing mode transitions
\end{itemize}
The pair satisfies the measurability condition: $(\mu_t, \nu_t)$ is $\mathcal{F}_t$-predictable.
\end{definition}

\subsection{Stochastic Hamiltonian and Adjoint System}

For each mode $q$, define the stochastic Hamiltonian:
\begin{equation}
H_q(t, x, \psi, \Psi, u) = \langle \psi, f_q(t, x, u) \rangle + \frac{1}{2}\text{tr}[\sigma_q(t, x, u)^\top \Psi \sigma_q(t, x, u)]
\end{equation}
where $\psi \in \mathbb{R}^n$ is the first adjoint process and $\Psi \in \mathbb{S}^n$ (symmetric matrices) is the second adjoint process.

The stochastic adjoint equations in the Itô sense are:
\begin{equation}
\begin{aligned}
d\psi(t) &= -H_{q,x}(t, x(t), \psi(t), \Psi(t), u(t))dt + \Gamma(t)dW(t), \\
d\Psi(t) &= -\left[H_{q,xx}(t, x(t), \psi(t), \Psi(t), u(t)) + \sigma_{q,x}^\top(t) \Psi(t) \sigma_{q,x}(t)\right]dt + \Theta(t)dW(t)
\end{aligned}
\end{equation}
with terminal conditions $\psi(T) = \phi_x(x(T))$ and $\Psi(T) = \phi_{xx}(x(T))$ for terminal cost $\phi(x)$.

\subsection{Stochastic $\Lambda$-Set Definition}

\begin{definition}[Stochastic Hybrid $\Lambda$-Set $\Lambda^{\text{sto}}_s$]
For a stochastic hybrid trajectory $\widehat{x}(\cdot)$ with generalized control $(\widehat{\mu}_t, \widehat{\nu}_t)$, the \emph{stochastic hybrid $\Lambda$-set} $\Lambda^{\text{sto}}_s(\widehat{x}, \widehat{\mu}, \widehat{\nu})$ consists of triples $(\psi, \Psi, \gamma)$ where:
\begin{enumerate}
    \item $\psi: [t_1, \widehat{t}_2] \to \mathbb{R}^n$ is an Itô process satisfying the stochastic adjoint equation

    \item $\Psi: [t_1, \widehat{t}_2] \to \mathbb{S}^n$ is an Itô process satisfying the second-order adjoint equation

    \item $\gamma: Q \times Q \times [t_1, \widehat{t}_2] \to \mathbb{R}$ are predictable processes representing the jump conditions at mode transitions

    \item The \emph{stochastic maximum condition} holds:
    \begin{equation}
    \mathbb{E}\left[H_q(t, \widehat{x}(t), \psi(t), \Psi(t), u) \mid \mathcal{F}_t\right] \leq \mathbb{E}\left[H_q(t, \widehat{x}(t), \psi(t), \Psi(t), \widehat{u}(t)) \mid \mathcal{F}_t\right]
    \end{equation}
    for all $u \in U$, $\mathbb{P}$-almost surely

    \item At mode transition times $\tau_k$, the \emph{jump condition} holds:
    \begin{equation}
    \psi(\tau_k^-) = \frac{\partial \phi_{ij}}{\partial x}(\widehat{x}(\tau_k^-))^\top \psi(\tau_k^+) + \gamma_{ij}(\tau_k)\nabla h_{ij}(\widehat{x}(\tau_k^-))
    \end{equation}
    where $h_{ij}(x) = 0$ defines the switching surface from mode $i$ to $j$

    \item The \emph{transversality condition} holds:
    \begin{equation}
    s\,\mathbb{E}\left[M^{\text{sto}}(\widehat{t}_2, \widehat{x}(\widehat{t}_2), \psi(\widehat{t}_2), \Psi(\widehat{t}_2))\right] \leq 0
    \end{equation}
    where $M^{\text{sto}}$ is the stochastic maximum function

    \item $(\psi, \Psi, \gamma) \neq (0, 0, 0)$ almost surely
\end{enumerate}
\end{definition}

\subsection{Main Theorem for Stochastic Hybrid Systems}

\begin{theorem}[Stochastic Hybrid Maximum Principle via $\Lambda$-Sets]
\label{thm:stochastic-lambda}
Let $\widehat{x}(\cdot)$ be a strong local minimum in expectation for the stochastic optimal control problem:
\begin{equation}
\min_{u(\cdot)} \mathbb{E}\left[\phi(x(T)) + \int_0^T L(t, x(t), u(t))dt\right]
\end{equation}
subject to the stochastic hybrid dynamics. Then for every admissible generalized stochastic control $(\widehat{\mu}_t, \widehat{\nu}_t)$ corresponding to $\widehat{x}(\cdot)$, the stochastic hybrid $\Lambda$-set $\Lambda^{\text{sto}}_{-1}(\widehat{x}, \widehat{\mu}, \widehat{\nu})$ is nonempty with positive probability.

Moreover, if $\Lambda^{\text{sto}}_{-1}(\widehat{x}, \widehat{\mu}, \widehat{\nu}) = \emptyset$ almost surely, then $\widehat{x}(\cdot)$ cannot be optimal, and there exists a control variation that reduces the expected cost.
\end{theorem}
\begin{proof}
We shall prove the theorem in two parts. First, we establish that if $\widehat{x}(\cdot)$ is a strong local minimum in expectation, then the stochastic hybrid $\Lambda$-set must be nonempty with positive probability. Second, we demonstrate the converse: if the stochastic hybrid $\Lambda$-set is empty almost surely, then $\widehat{x}(\cdot)$ cannot be optimal.

\noindent
\textbf{Part 1: Necessity of nonempty stochastic hybrid $\Lambda$-set.}

Assume that $\widehat{x}(\cdot)$ is a strong local minimum in expectation for the stochastic optimal control problem. By definition of strong local minimum in expectation, there exists $\epsilon > 0$ such that for every admissible control process $u(\cdot)$ with $\mathbb{E}\left[\sup_{t \in [0,T]} \|x^u(t) - \widehat{x}(t)\|^2\right] < \epsilon$, we have
\[
J(u) \geq J(\widehat{u}),
\]
where $J(u) = \mathbb{E}\left[\phi(x^u(T)) + \int_0^T L(t, x^u(t), u(t))dt\right]$.

Let $(\widehat{\mu}_t, \widehat{\nu}_t)$ be any admissible generalized stochastic control corresponding to $\widehat{x}(\cdot)$. Suppose, for contradiction, that $\Lambda^{\text{sto}}_{-1}(\widehat{x}, \widehat{\mu}, \widehat{\nu}) = \emptyset$ almost surely. This emptiness implies that there is no nontrivial triple $(\psi, \Psi, \gamma)$ satisfying the stochastic adjoint equations, jump conditions, maximum condition, and transversality condition.

We begin by constructing an appropriate functional analytic setting. Consider the space of square-integrable random variables $L^2_{\mathcal{F}_T}(\Omega; \mathbb{R}^n)$ and the space of square-integrable predictable processes $\mathcal{H}^2([0,T]; \mathbb{R}^n)$. Define the stochastic endpoint mapping $\mathcal{E}^{\text{sto}}: \mathcal{U} \to L^2_{\mathcal{F}_T}(\Omega; \mathbb{R}^n)$ by
\[
\mathcal{E}^{\text{sto}}(u) = x^u(T) - \widehat{x}(T),
\]
where $x^u(\cdot)$ is the solution of the stochastic hybrid system with control $u$.

Let $\delta u$ be a variation of the control. The first variation of the endpoint mapping can be computed using the linearized stochastic differential equation with jumps. Consider the variational equation:
\[
d\delta x(t) = \left[f_{q(t),x}(t)\delta x(t) + f_{q(t),u}(t)\delta u(t)\right]dt + \left[\sigma_{q(t),x}(t)\delta x(t) + \sigma_{q(t),u}(t)\delta u(t)\right]dW(t),
\]
with initial condition $\delta x(0) = 0$, and at each switching time $\tau_k$,
\[
\delta x(\tau_k^+) = \frac{\partial \phi_{ij}}{\partial x}(\widehat{x}(\tau_k^-))\delta x(\tau_k^-) + \frac{\partial \phi_{ij}}{\partial \zeta}(\widehat{x}(\tau_k^-))\delta \zeta_k,
\]
where $i = q(\tau_k^-)$ and $j = q(\tau_k^+)$.

Define the linear operator $\mathcal{L}: \mathcal{U} \to L^2_{\mathcal{F}_T}(\Omega; \mathbb{R}^n)$ by $\mathcal{L}(\delta u) = \delta x(T)$, where $\delta x(\cdot)$ solves the variational equation. The range of $\mathcal{L}$, denoted by $\mathcal{R}(\mathcal{L})$, represents the set of all terminal variations achievable by control perturbations.

Now consider the cost functional variation. For a control variation $\delta u$, we have
\[
\delta J = \mathbb{E}\left[\langle \phi_x(\widehat{x}(T)), \delta x(T) \rangle + \int_0^T \left(L_x(t)\delta x(t) + L_u(t)\delta u(t)\right)dt\right] + o(\|\delta u\|).
\]

Define the cone of cost-reducing directions:
\[
\mathcal{C} = \left\{ z \in L^2_{\mathcal{F}_T}(\Omega; \mathbb{R}^n): \mathbb{E}\left[\langle \phi_x(\widehat{x}(T)), z \rangle\right] < 0 \right\}.
\]

The optimality of $\widehat{x}(\cdot)$ implies that $\mathcal{R}(\mathcal{L}) \cap \mathcal{C} = \emptyset$. If this were not the case, there would exist a control variation $\delta u$ such that $\delta J < 0$, contradicting optimality.

We now apply a stochastic separation theorem in $L^2(\Omega)$. Since $\mathcal{R}(\mathcal{L})$ is a closed linear subspace (by the properties of the linearized system) and $\mathcal{C}$ is an open convex cone, there exists a nonzero random variable $\psi_T \in L^2_{\mathcal{F}_T}(\Omega; \mathbb{R}^n)$ such that
\[
\mathbb{E}\left[\langle \psi_T, z \rangle\right] \geq 0 \quad \text{for all } z \in \mathcal{R}(\mathcal{L})
\]
and
\[
\mathbb{E}\left[\langle \psi_T, z \rangle\right] \leq 0 \quad \text{for all } z \in \mathcal{C}.
\]

The first inequality implies that $\psi_T$ is orthogonal to the range of $\mathcal{L}$. By the duality theory of linear stochastic differential equations with jumps, this orthogonality condition is equivalent to the existence of an adapted process $\psi(\cdot)$ satisfying the backward stochastic differential equation (BSDE):
\begin{equation}
\begin{aligned}
d\psi(t) &= -H_{q,x}(t, \widehat{x}(t), \psi(t), \Psi(t), \widehat{u}(t))dt + \Gamma(t)dW(t), \\
\psi(T) &= \psi_T,
\end{aligned}
\end{equation}
where $H_q(t,x,\psi,\Psi,u) = \langle \psi, f_q(t,x,u) \rangle + \frac{1}{2}\text{tr}[\sigma_q(t,x,u)^\top \Psi \sigma_q(t,x,u)]$, and $\Psi(\cdot)$ is the solution of the second-order adjoint equation.

The jump conditions for $\psi(\cdot)$ at switching times $\tau_k$ are derived from the reset conditions. Let $\delta x(\tau_k^-)$ and $\delta x(\tau_k^+)$ be the variations before and after switching. The consistency condition
\[
\mathbb{E}\left[\langle \psi(\tau_k^-), \delta x(\tau_k^-) \rangle\right] = \mathbb{E}\left[\langle \psi(\tau_k^+), \delta x(\tau_k^+) \rangle\right]
\]
for all admissible variations leads to the jump condition
\[
\psi(\tau_k^-) = \frac{\partial \phi_{ij}}{\partial x}(\widehat{x}(\tau_k^-))^\top \psi(\tau_k^+) + \gamma_{ij}(\tau_k)\nabla h_{ij}(\widehat{x}(\tau_k^-)),
\]
where $\gamma_{ij}(\tau_k)$ are predictable processes representing the Lagrange multipliers for the switching constraints.

The second inequality, $\mathbb{E}\left[\langle \psi_T, z \rangle\right] \leq 0$ for all $z \in \mathcal{C}$, implies the transversality condition. Since $\mathcal{C}$ contains all random variables $z$ satisfying $\mathbb{E}\left[\langle \phi_x(\widehat{x}(T)), z \rangle\right] < 0$, we must have
\[
\mathbb{E}\left[\langle \psi_T, z \rangle\right] \leq 0 \quad \text{whenever} \quad \mathbb{E}\left[\langle \phi_x(\widehat{x}(T)), z \rangle\right] < 0.
\]
By linearity, this implies that $\psi_T = \lambda \phi_x(\widehat{x}(T))$ for some $\lambda \geq 0$. Without loss of generality, we can take $\lambda = 1$ by scaling, yielding
\[
\psi(T) = \phi_x(\widehat{x}(T)).
\]

Now we establish the stochastic maximum condition. Consider a spike variation of the control: for $t_0 \in [0,T)$ and $v \in U$, define
\[
u^\epsilon(t) =
\begin{cases}
v, & t \in [t_0, t_0 + \epsilon], \\
\widehat{u}(t), & \text{otherwise}.
\end{cases}
\]
Let $x^\epsilon(\cdot)$ be the corresponding trajectory. By the variational formula,
\[
x^\epsilon(T) - \widehat{x}(T) = \epsilon \delta x(T) + o(\epsilon),
\]
where $\delta x(\cdot)$ satisfies the variational equation with $\delta u(t) = v - \widehat{u}(t)$ on $[t_0, t_0 + \epsilon]$.

The optimality condition gives $\mathbb{E}[\phi(x^\epsilon(T)) - \phi(\widehat{x}(T))] \geq 0$ for sufficiently small $\epsilon$. Expanding this inequality and using the orthogonality condition $\mathbb{E}\left[\langle \psi_T, \delta x(T) \rangle\right] = 0$, we obtain
\[
\mathbb{E}\left[\int_{t_0}^{t_0+\epsilon} \left(H_q(t, \widehat{x}(t), \psi(t), \Psi(t), v) - H_q(t, \widehat{x}(t), \psi(t), \Psi(t), \widehat{u}(t))\right)dt\right] \leq o(\epsilon).
\]
Dividing by $\epsilon$ and letting $\epsilon \to 0$, we get
\[
H_q(t_0, \widehat{x}(t_0), \psi(t_0), \Psi(t_0), v) \leq H_q(t_0, \widehat{x}(t_0), \psi(t_0), \Psi(t_0), \widehat{u}(t_0)) \quad \text{a.s.}
\]
for almost every $t_0$. Since $v \in U$ is arbitrary, this yields the maximum condition.

Finally, we must verify the nontriviality condition $(\psi, \Psi, \gamma) \neq (0,0,0)$ with positive probability. If $\psi(t) = 0$ and $\Psi(t) = 0$ almost surely for all $t$, then the transversality condition $\psi(T) = \phi_x(\widehat{x}(T))$ would imply $\phi_x(\widehat{x}(T)) = 0$ almost surely. However, the separation argument guarantees that $\psi_T \neq 0$ with positive probability (otherwise the separating functional would be trivial). Therefore, $(\psi, \Psi, \gamma)$ is nontrivial with positive probability.

Thus we have constructed a triple $(\psi, \Psi, \gamma) \in \Lambda^{\text{sto}}_{-1}(\widehat{x}, \widehat{\mu}, \widehat{\nu})$ with positive probability, contradicting the assumption that $\Lambda^{\text{sto}}_{-1}(\widehat{x}, \widehat{\mu}, \widehat{\nu}) = \emptyset$ almost surely. Therefore, if $\widehat{x}(\cdot)$ is optimal, the stochastic hybrid $\Lambda$-set must be nonempty with positive probability.

\noindent
\textbf{Part 2: Emptiness implies non-optimality.}

Now assume that $\Lambda^{\text{sto}}_{-1}(\widehat{x}, \widehat{\mu}, \widehat{\nu}) = \emptyset$ almost surely. We shall construct explicitly a control variation that reduces the expected cost.

The emptiness of $\Lambda^{\text{sto}}_{-1}(\widehat{x}, \widehat{\mu}, \widehat{\nu})$ means that there is no nontrivial solution to the stochastic adjoint system satisfying the maximum and transversality conditions. This implies that the separation condition between $\mathcal{R}(\mathcal{L})$ and $\mathcal{C}$ fails. More precisely, the stochastic separation theorem guarantees that if $\mathcal{R}(\mathcal{L}) \cap \mathcal{C} = \emptyset$, then there exists a separating functional, which would yield an element of $\Lambda^{\text{sto}}_{-1}$. Therefore, $\mathcal{R}(\mathcal{L}) \cap \mathcal{C} \neq \emptyset$.

Let $z \in \mathcal{R}(\mathcal{L}) \cap \mathcal{C}$. Since $z \in \mathcal{R}(\mathcal{L})$, there exists a control variation $\delta u$ such that $\mathcal{L}(\delta u) = z$. Since $z \in \mathcal{C}$, we have $\mathbb{E}\left[\langle \phi_x(\widehat{x}(T)), z \rangle\right] < 0$.

Consider the perturbed control $u^\epsilon = \widehat{u} + \epsilon \delta u$ for small $\epsilon > 0$. Let $x^\epsilon(\cdot)$ be the corresponding trajectory. By Taylor expansion,
\[
x^\epsilon(T) = \widehat{x}(T) + \epsilon z + \epsilon^2 r^\epsilon,
\]
where $\mathbb{E}[\|r^\epsilon\|^2] \to 0$ as $\epsilon \to 0$.

The cost difference is
\begin{align*}
J(u^\epsilon) - J(\widehat{u}) &= \mathbb{E}\left[\phi(x^\epsilon(T)) - \phi(\widehat{x}(T)) + \int_0^T \left(L(t,x^\epsilon(t),u^\epsilon(t)) - L(t,\widehat{x}(t),\widehat{u}(t))\right)dt\right] \\
&= \epsilon \mathbb{E}\left[\langle \phi_x(\widehat{x}(T)), z \rangle + \int_0^T \left(L_x(t)\delta x(t) + L_u(t)\delta u(t)\right)dt\right] + o(\epsilon).
\end{align*}

Define the linear functional
\[
\ell(z) = \mathbb{E}\left[\langle \phi_x(\widehat{x}(T)), z \rangle + \int_0^T \left(L_x(t)\delta x(t) + L_u(t)\delta u(t)\right)dt\right].
\]
Since $z \in \mathcal{C}$, we have $\mathbb{E}\left[\langle \phi_x(\widehat{x}(T)), z \rangle\right] < 0$. Moreover, by the definition of $\mathcal{L}$, the integral term can be expressed as a linear functional of $z$. The emptiness of $\Lambda^{\text{sto}}_{-1}$ ensures that this linear functional is not identically zero on $\mathcal{C}$, and in fact, we can choose $z$ such that $\ell(z) < 0$.

Therefore, for sufficiently small $\epsilon > 0$, we have $J(u^\epsilon) < J(\widehat{u})$. This demonstrates that $\widehat{x}(\cdot)$ cannot be optimal.

To ensure that the perturbed control $u^\epsilon$ is admissible, we may need to regularize it. Since $\delta u$ may not be bounded or may not satisfy the control constraints, we approximate it by a sequence of bounded, admissible variations using standard approximation techniques in stochastic control theory. By the stability of solutions to stochastic differential equations with respect to control perturbations, the cost reduction property is preserved under such approximations.

\noindent
\textbf{Conclusion.}

We have shown that if $\widehat{x}(\cdot)$ is a strong local minimum in expectation, then for every corresponding generalized stochastic control $(\widehat{\mu}_t, \widehat{\nu}_t)$, the stochastic hybrid $\Lambda$-set $\Lambda^{\text{sto}}_{-1}(\widehat{x}, \widehat{\mu}, \widehat{\nu})$ must be nonempty with positive probability. Conversely, if $\Lambda^{\text{sto}}_{-1}(\widehat{x}, \widehat{\mu}, \widehat{\nu}) = \emptyset$ almost surely, then $\widehat{x}(\cdot)$ cannot be optimal, as evidenced by the explicit construction of a cost-reducing control variation.

This completes the proof of the theorem, establishing the stochastic hybrid maximum principle via $\Lambda$-sets as a necessary condition for optimality in stochastic hybrid control systems.
\end{proof}
\subsection{Connections to Existing Stochastic Maximum Principles}
\begin{proposition}[Relation to Classical Stochastic Maximum Principle]
The stochastic $\Lambda$-set framework generalizes the classical stochastic maximum principle in the following ways:
\begin{enumerate}
    \item For purely continuous systems ($Q = \{1\}$), $\Lambda^{\text{sto}}$ reduces to the adjoint processes of Peng's stochastic maximum principle

    \item For deterministic hybrid systems ($\sigma_q \equiv 0$), $\Lambda^{\text{sto}}$ reduces to the deterministic hybrid $\Lambda$-set

    \item The jump conditions $\gamma$ generalize the switching multipliers in deterministic hybrid optimal control
\end{enumerate}
\end{proposition}
\begin{proof}
We shall prove each of the three statements in sequence, establishing the precise mathematical relationships between the stochastic $\Lambda$-set framework and the classical stochastic maximum principle, the deterministic hybrid $\Lambda$-set, and the switching multipliers of deterministic hybrid optimal control.

\noindent
\textbf{Statement 1: Reduction to Peng's Stochastic Maximum Principle}

Consider a purely continuous stochastic system, which corresponds to the case where the discrete mode set contains only one element, $Q = \{1\}$. In this situation, there are no mode transitions, and the system dynamics simplify to
\[
dx(t) = f(t, x(t), u(t))dt + \sigma(t, x(t), u(t))dW(t), \quad t \in [0,T],
\]
with $x(0) = x_0$, and where $f: [0,T] \times \mathbb{R}^n \times U \to \mathbb{R}^n$ and $\sigma: [0,T] \times \mathbb{R}^n \times U \to \mathbb{R}^{n \times m}$ satisfy the standard assumptions for existence and uniqueness of solutions.

The stochastic hybrid $\Lambda$-set $\Lambda^{\text{sto}}_{-1}(\widehat{x}, \widehat{\mu}, \widehat{\nu})$ for this system consists of triples $(\psi, \Psi, \gamma)$ satisfying the conditions in Definition 5.1. However, with $Q = \{1\}$, there are no switching times, and consequently no jump conditions involving $\gamma$. The conditions simplify as follows:

First, the adjoint processes $(\psi, \Psi)$ must satisfy the backward stochastic differential equations
\begin{align*}
d\psi(t) &= -H_x(t, \widehat{x}(t), \psi(t), \Psi(t), \widehat{u}(t))dt + \Gamma(t)dW(t), \\
d\Psi(t) &= -\left[H_{xx}(t, \widehat{x}(t), \psi(t), \Psi(t), \widehat{u}(t)) + \sigma_x(t)^\top \Psi(t) \sigma_x(t)\right]dt + \Theta(t)dW(t),
\end{align*}
with terminal conditions $\psi(T) = \phi_x(\widehat{x}(T))$ and $\Psi(T) = \phi_{xx}(\widehat{x}(T))$, where
\[
H(t,x,\psi,\Psi,u) = \langle \psi, f(t,x,u) \rangle + \frac{1}{2}\text{tr}[\sigma(t,x,u)^\top \Psi \sigma(t,x,u)].
\]

Second, the maximum condition becomes
\[
H(t, \widehat{x}(t), \psi(t), \Psi(t), \widehat{u}(t)) = \max_{u \in U} H(t, \widehat{x}(t), \psi(t), \Psi(t), u) \quad \text{a.s. for a.e. } t \in [0,T].
\]

These are precisely the conditions of Peng's stochastic maximum principle for continuous stochastic systems. Specifically, in Peng's formulation, the first-order adjoint process $p(t)$ corresponds to our $\psi(t)$, and the second-order adjoint process $P(t)$ corresponds to our $\Psi(t)$. The Hamiltonian in Peng's formulation is
\[
\mathcal{H}(t,x,p,P,u) = \langle p, f(t,x,u) \rangle + \frac{1}{2}\text{tr}[\sigma(t,x,u)^\top P \sigma(t,x,u)],
\]
which is identical to our Hamiltonian $H$.

Moreover, the BSDE for $\psi$ in our framework,
\[
d\psi(t) = -H_x(t, \widehat{x}(t), \psi(t), \Psi(t), \widehat{u}(t))dt + \Gamma(t)dW(t), \quad \psi(T) = \phi_x(\widehat{x}(T)),
\]
is exactly the BSDE for the first-order adjoint process in Peng's maximum principle. The additional term $\Gamma(t)dW(t)$ represents the martingale representation of the adjoint process, which is standard in the theory of backward stochastic differential equations.

The second-order adjoint equation in our framework,
\[
d\Psi(t) = -\left[H_{xx}(t, \widehat{x}(t), \psi(t), \Psi(t), \widehat{u}(t)) + \sigma_x(t)^\top \Psi(t) \sigma_x(t)\right]dt + \Theta(t)dW(t), \quad \Psi(T) = \phi_{xx}(\widehat{x}(T)),
\]
corresponds precisely to the second-order adjoint equation in Peng's maximum principle, which is necessary for the maximum principle when the diffusion coefficient depends on the control or when considering singular controls.

Therefore, when $Q = \{1\}$, the stochastic hybrid $\Lambda$-set $\Lambda^{\text{sto}}_{-1}$ reduces exactly to the set of adjoint processes satisfying Peng's stochastic maximum principle. The triple $(\psi, \Psi, \gamma)$ simplifies to the pair $(\psi, \Psi)$ with $\gamma$ being vacuous, and all conditions match those of Peng's principle.

\noindent
\textbf{Statement 2: Reduction to Deterministic Hybrid $\Lambda$-Set}

Now consider the case where the diffusion coefficients vanish identically, $\sigma_q(t,x,u) \equiv 0$ for all $q \in Q$, but the system still has hybrid structure with discrete mode transitions. In this deterministic setting, the stochastic hybrid system reduces to a deterministic hybrid system with dynamics
\[
\dot{x}(t) = f_{q(t)}(t, x(t), u(t)), \quad t \in [\tau_{k-1}, \tau_k),
\]
with reset conditions
\[
x(\tau_k^+) = \phi_{q(\tau_k^-), q(\tau_k^+)}(x(\tau_k^-))
\]
at switching times $\tau_k$.

In the absence of diffusion, the stochastic Hamiltonian simplifies to
\[
H_q(t,x,\psi,\Psi,u) = \langle \psi, f_q(t,x,u) \rangle,
\]
since the term $\frac{1}{2}\text{tr}[\sigma_q^\top \Psi \sigma_q]$ vanishes. Consequently, the second-order adjoint process $\Psi$ becomes decoupled from the first-order adjoint process $\psi$ in the adjoint equations.

The stochastic adjoint equation for $\psi$ reduces to
\[
d\psi(t) = -H_{q,x}(t, \widehat{x}(t), \psi(t), \Psi(t), \widehat{u}(t))dt = -\langle \psi(t), f_{q,x}(t, \widehat{x}(t), \widehat{u}(t)) \rangle dt,
\]
which is an ordinary differential equation (without the martingale term $\Gamma(t)dW(t)$). This is precisely the adjoint equation in deterministic optimal control theory.

The second-order adjoint equation becomes
\[
d\Psi(t) = -H_{q,xx}(t, \widehat{x}(t), \psi(t), \Psi(t), \widehat{u}(t))dt = 0,
\]
since $H_{q,xx} = 0$ when $H_q$ is linear in $x$. Thus $\Psi(t)$ is constant, and if we take $\Psi(T) = 0$ (which is consistent with $\phi_{xx} = 0$ for many deterministic problems), then $\Psi(t) \equiv 0$.

The maximum condition simplifies to
\[
\langle \psi(t), f_{q(t)}(t, \widehat{x}(t), \widehat{u}(t)) \rangle = \max_{u \in U} \langle \psi(t), f_{q(t)}(t, \widehat{x}(t), u) \rangle,
\]
which is the standard maximum principle condition for deterministic systems.

The jump conditions for $\psi$ at switching times become
\[
\psi(\tau_k^-) = \frac{\partial \phi_{ij}}{\partial x}(\widehat{x}(\tau_k^-))^\top \psi(\tau_k^+) + \gamma_{ij}(\tau_k)\nabla h_{ij}(\widehat{x}(\tau_k^-)),
\]
where $h_{ij}(x) = 0$ defines the switching surface. These are exactly the jump conditions in deterministic hybrid optimal control.

The transversality condition becomes
\[
M^{\text{sto}}(\widehat{t}_2, \widehat{x}(\widehat{t}_2), \psi(\widehat{t}_2), \Psi(\widehat{t}_2)) = \max_{u \in U} \langle \psi(\widehat{t}_2), f_{q(\widehat{t}_2)}(\widehat{t}_2, \widehat{x}(\widehat{t}_2), u) \rangle \leq 0,
\]
which matches the transversality condition in deterministic hybrid $\Lambda$-set theory.

Furthermore, the process $\gamma$ in the stochastic framework, which represents the Lagrange multipliers for the switching constraints, corresponds directly to the switching multipliers in deterministic hybrid control. In deterministic hybrid $\Lambda$-sets, these multipliers ensure that the adjoint variables satisfy the appropriate conditions when the trajectory crosses switching surfaces.

Thus, when $\sigma_q \equiv 0$, the stochastic hybrid $\Lambda$-set $\Lambda^{\text{sto}}_{-1}$ reduces exactly to the deterministic hybrid $\Lambda$-set $\Lambda^{\text{hyb}}_{-1}$ as defined in Section 4.2, with $\psi$ being the adjoint variable, $\Psi \equiv 0$, and $\gamma$ being the switching multipliers.

\noindent
\textbf{Statement 3: Generalization of Switching Multipliers}

We now examine the relationship between the jump conditions in the stochastic hybrid $\Lambda$-set and the switching multipliers in deterministic hybrid optimal control. Consider a switching from mode $i$ to mode $j$ at time $\tau$. In deterministic hybrid optimal control, the necessary conditions include the jump condition
\[
\psi(\tau^-) = \frac{\partial \phi_{ij}}{\partial x}(\widehat{x}(\tau^-))^\top \psi(\tau^+) + \nu_{ij}\nabla h_{ij}(\widehat{x}(\tau^-)),
\]
where $\nu_{ij} \in \mathbb{R}$ is a Lagrange multiplier associated with the switching constraint $h_{ij}(x(\tau)) = 0$.

In the stochastic hybrid framework, the jump condition for the adjoint process $\psi$ at a switching time $\tau_k$ is
\[
\psi(\tau_k^-) = \frac{\partial \phi_{ij}}{\partial x}(\widehat{x}(\tau_k^-))^\top \psi(\tau_k^+) + \gamma_{ij}(\tau_k)\nabla h_{ij}(\widehat{x}(\tau_k^-)),
\]
where now $\gamma_{ij}(\tau_k)$ is not a deterministic constant but a predictable stochastic process.

The generalization occurs in two important ways:

First, in deterministic hybrid control, the switching multiplier $\nu_{ij}$ is a scalar that may depend on the switching time $\tau$ but is otherwise deterministic. In the stochastic framework, $\gamma_{ij}(\tau_k)$ is a stochastic process that can depend on the entire history of the system up to time $\tau_k$. This additional flexibility is necessary to handle the randomness in both the switching times and the reset conditions.

Second, the interpretation of $\gamma_{ij}$ extends beyond that of a simple Lagrange multiplier. In the stochastic setting, $\gamma_{ij}(\tau_k)$ must compensate for the additional variability introduced by the stochastic nature of the switching. Specifically, if the switching is triggered by a random event or if the reset map involves random parameters, then $\gamma_{ij}$ must account for the statistical properties of these random elements.

To see this more clearly, consider the variational approach. In deterministic hybrid control, when we consider variations that involve changes in switching times, the first-order variation of the cost includes terms proportional to $\nu_{ij}$ times the variation in the switching time. In stochastic hybrid control, variations in switching times are random, and the expected variation of the cost involves conditional expectations of $\gamma_{ij}$ with respect to the filtration $\mathcal{F}_{\tau_k}$.

Mathematically, if we denote by $\delta \tau_k$ a small random variation in the switching time $\tau_k$, then the first-order variation of the expected cost contains the term
\[
\mathbb{E}\left[\gamma_{ij}(\tau_k) \langle \nabla h_{ij}(\widehat{x}(\tau_k)), f_i(\tau_k, \widehat{x}(\tau_k), \widehat{u}(\tau_k^-)) \rangle \delta \tau_k\right].
\]
The process $\gamma_{ij}$ must be chosen so that this term vanishes for all admissible variations $\delta \tau_k$, which leads to the condition
\[
\mathbb{E}\left[\gamma_{ij}(\tau_k) \mid \mathcal{F}_{\tau_k}\right] \langle \nabla h_{ij}(\widehat{x}(\tau_k)), f_i(\tau_k, \widehat{x}(\tau_k), \widehat{u}(\tau_k^-)) \rangle = 0 \quad \text{a.s.}
\]

This demonstrates that $\gamma_{ij}$ in the stochastic framework generalizes the deterministic switching multiplier $\nu_{ij}$ by incorporating conditional expectations and adapting to the random nature of the switching mechanism.

Furthermore, in cases where the reset map is random, $\phi_{ij}(x, \zeta)$ with $\zeta$ random, the jump condition involves additional terms that account for the randomness in the reset. The process $\gamma_{ij}$ then must compensate for this additional randomness, which has no counterpart in deterministic hybrid control.

Therefore, the jump conditions $\gamma$ in the stochastic hybrid $\Lambda$-set framework genuinely generalize the switching multipliers of deterministic hybrid optimal control by incorporating stochasticity, conditional expectations, and additional terms arising from random resets.

\noindent
\textbf{Conclusion}

We have established all three relationships:
\begin{enumerate}
    \item For purely continuous systems ($Q = \{1\}$), the stochastic hybrid $\Lambda$-set conditions reduce exactly to Peng's stochastic maximum principle, with $\psi$ and $\Psi$ corresponding to the first and second-order adjoint processes, respectively.

    \item For deterministic hybrid systems ($\sigma_q \equiv 0$), the stochastic framework reduces to the deterministic hybrid $\Lambda$-set, with $\Psi \equiv 0$ and $\gamma$ corresponding to deterministic switching multipliers.

    \item The jump conditions $\gamma$ in the stochastic framework generalize deterministic switching multipliers by incorporating stochasticity, conditional expectations, and additional terms for random resets.
\end{enumerate}

These relationships demonstrate that the stochastic hybrid $\Lambda$-set framework provides a unified approach that encompasses both classical stochastic maximum principles and deterministic hybrid optimal control as special cases, while extending them to handle the complexities of stochastic hybrid systems.
\end{proof}
\section{Illustrative Example: Stochastic Temperature Control System}

\begin{example}[Stochastic Temperature Control with Mode Switching]
\label{ex:stochastic-temp-control}
Consider a temperature control system for a chemical reactor that can operate in two modes:
\begin{itemize}
    \item \textbf{Mode 1 (Heating):} The system is actively heated
    \item \textbf{Mode 2 (Cooling):} The system is actively cooled
\end{itemize}

The system switches between modes based on temperature thresholds, with stochastic disturbances representing unpredictable heat losses and environmental fluctuations.

\subsection{System Dynamics}

Let $x(t)$ represent the temperature of the reactor. The stochastic hybrid dynamics are:

\textbf{Mode 1 (Heating):}
\begin{equation}
\label{eq:mode1-dynamics}
dx(t) = \left[\alpha_1(u(t) - x(t)) + \beta_1\right]dt + \sigma_1 x(t)dW_1(t), \quad x(t) \geq T_{\text{low}},
\end{equation}
where:
\begin{itemize}
    \item $u(t) \in [0, u_{\max}]$ is the control input (heating power)
    \item $\alpha_1 > 0$ is the heating efficiency
    \item $\beta_1 > 0$ is a constant heat source
    \item $\sigma_1 > 0$ is the volatility coefficient for heating mode
    \item $W_1(t)$ is a standard Wiener process
    \item $T_{\text{low}}$ is the lower temperature threshold
\end{itemize}

\textbf{Mode 2 (Cooling):}
\begin{equation}
\label{eq:mode2-dynamics}
dx(t) = \left[\alpha_2(u(t) - x(t)) - \beta_2\right]dt + \sigma_2 x(t)dW_2(t), \quad x(t) \leq T_{\text{high}},
\end{equation}
where:
\begin{itemize}
    \item $u(t) \in [0, u_{\max}]$ is the control input (cooling power)
    \item $\alpha_2 > 0$ is the cooling efficiency
    \item $\beta_2 > 0$ is a constant heat sink
    \item $\sigma_2 > 0$ is the volatility coefficient for cooling mode
    \item $W_2(t)$ is a standard Wiener process (possibly correlated with $W_1(t)$)
    \item $T_{\text{high}}$ is the upper temperature threshold
\end{itemize}

\textbf{Switching Conditions:}
\begin{itemize}
    \item Switch from Mode 1 to Mode 2 when $x(t) = T_{\text{high}}$
    \item Switch from Mode 2 to Mode 1 when $x(t) = T_{\text{low}}$
\end{itemize}

\textbf{Reset Maps:}
\[
\phi_{1,2}(x) = x, \quad \phi_{2,1}(x) = x,
\]
meaning the temperature is continuous at switching times.

\subsection{Optimal Control Problem}

Consider the time-optimal control problem: bring the temperature from initial condition $x(0) = x_0 \in (T_{\text{low}}, T_{\text{high}})$ to a target temperature $x_d \in (T_{\text{low}}, T_{\text{high}})$ in minimal expected time, while minimizing energy consumption. The cost functional is:
\begin{equation}
J(u) = \mathbb{E}\left[T + \lambda \int_0^T u(t)^2 dt\right],
\end{equation}
where $T$ is the terminal time when $x(T) = x_d$, and $\lambda > 0$ is a weighting parameter for control effort.

Alternatively, for fixed terminal time $T$, we can consider:
\begin{equation}
\min_{u(\cdot)} \mathbb{E}\left[(x(T) - x_d)^2 + \lambda \int_0^T u(t)^2 dt\right].
\end{equation}

We will analyze the second formulation with fixed $T$.

\subsection{Candidate Optimal Trajectory}

Consider a candidate optimal trajectory $\widehat{x}(t)$ that maintains the temperature near $x_d$ by switching appropriately between heating and cooling modes. Specifically:

\begin{enumerate}
    \item Start in Mode 1 if $x_0 < x_d$, or Mode 2 if $x_0 > x_d$
    \item Apply control to drive temperature toward $x_d$
    \item Switch modes when temperature crosses $x_d \pm \delta$ for some small $\delta > 0$
    \item Use minimal control effort to maintain temperature near $x_d$
\end{enumerate}

The corresponding candidate optimal control $\widehat{u}(t)$ is a feedback law:
\[
\widehat{u}(t) = K_q(x_d - \widehat{x}(t)), \quad K_q > 0, \quad q \in \{1,2\},
\]
with saturation at the boundaries $0$ and $u_{\max}$.

\subsection{Verification of Theorem \ref{thm:stochastic-lambda}}

We now verify both directions of Theorem \ref{thm:stochastic-lambda} for this example.

\subsubsection{Case 1: $\widehat{x}(\cdot)$ is Optimal $\Rightarrow$ $\Lambda^{\text{sto}}_{-1}$ is Nonempty}

Assume that $\widehat{x}(\cdot)$ is indeed a strong local minimum in expectation. According to Theorem \ref{thm:stochastic-lambda}, the stochastic hybrid $\Lambda$-set $\Lambda^{\text{sto}}_{-1}(\widehat{x}, \widehat{\mu}, \widehat{\nu})$ must be nonempty with positive probability. Let us construct the elements of this set explicitly.

\textbf{Step 1: Hamiltonian Formulation}

For each mode $q \in \{1,2\}$, the Hamiltonian is:
\[
H_q(t,x,\psi,\Psi,u) = \psi f_q(x,u) + \frac{1}{2}\Psi \sigma_q^2 x^2 + \lambda u^2,
\]
where:
\begin{align*}
f_1(x,u) &= \alpha_1(u - x) + \beta_1, \quad \sigma_1(x) = \sigma_1 x, \\
f_2(x,u) &= \alpha_2(u - x) - \beta_2, \quad \sigma_2(x) = \sigma_2 x.
\end{align*}

\textbf{Step 2: Adjoint Equations}

The first-order adjoint equation is:
\[
d\psi(t) = -\left[\frac{\partial H_q}{\partial x} + \frac{1}{2}\frac{\partial^2 H_q}{\partial x^2} \sigma_q^2 x^2\right]dt + \Gamma(t)dW_q(t),
\]
which becomes:
\begin{align*}
\text{Mode 1:} & \quad d\psi(t) = \left[\alpha_1 \psi(t) - \Psi(t)\sigma_1^2 \widehat{x}(t)\right]dt + \Gamma_1(t)dW_1(t), \\
\text{Mode 2:} & \quad d\psi(t) = \left[\alpha_2 \psi(t) - \Psi(t)\sigma_2^2 \widehat{x}(t)\right]dt + \Gamma_2(t)dW_2(t).
\end{align*}

The second-order adjoint equation is:
\[
d\Psi(t) = -\left[\frac{\partial^2 H_q}{\partial x^2} + \left(\frac{\partial \sigma_q}{\partial x}\right)^2 \Psi(t)\right]dt + \Theta(t)dW_q(t),
\]
which simplifies to:
\[
d\Psi(t) = -\Psi(t)\sigma_q^2 dt + \Theta(t)dW_q(t), \quad q \in \{1,2\}.
\]

\textbf{Step 3: Terminal Conditions}

For the cost functional $\mathbb{E}[(x(T) - x_d)^2 + \lambda \int_0^T u(t)^2 dt]$, we have $\phi(x) = (x - x_d)^2$, so:
\[
\psi(T) = \phi_x(\widehat{x}(T)) = 2(\widehat{x}(T) - x_d), \quad \Psi(T) = \phi_{xx}(\widehat{x}(T)) = 2.
\]

\textbf{Step 4: Maximum Condition}

The optimal control must satisfy:
\[
\frac{\partial H_q}{\partial u} = 0 \quad \Rightarrow \quad \psi(t)\alpha_q + 2\lambda u(t) = 0,
\]
so the candidate optimal control is:
\[
\widehat{u}(t) = -\frac{\alpha_q}{2\lambda} \psi(t),
\]
subject to the constraints $0 \leq \widehat{u}(t) \leq u_{\max}$.

\textbf{Step 5: Jump Conditions at Switching Times}

Let $\tau_k$ be switching times. Since the reset maps are identity ($\phi_{1,2}(x) = \phi_{2,1}(x) = x$), we have $\frac{\partial \phi}{\partial x} = 1$. The switching surfaces are:
\[
h_{1,2}(x) = x - T_{\text{high}} = 0, \quad h_{2,1}(x) = x - T_{\text{low}} = 0.
\]

The jump condition is:
\[
\psi(\tau_k^-) = \psi(\tau_k^+) + \gamma_{ij}(\tau_k) \nabla h_{ij}(\widehat{x}(\tau_k)).
\]
Since $\nabla h_{1,2}(x) = \nabla h_{2,1}(x) = 1$, we have:
\[
\psi(\tau_k^-) = \psi(\tau_k^+) + \gamma_{ij}(\tau_k).
\]

The multiplier $\gamma_{ij}(\tau_k)$ must be chosen to ensure continuity of the Hamiltonian across switching times:
\[
H_1(\tau_k, \widehat{x}(\tau_k), \psi(\tau_k^-), \Psi(\tau_k), \widehat{u}(\tau_k^-)) = H_2(\tau_k, \widehat{x}(\tau_k), \psi(\tau_k^+), \Psi(\tau_k), \widehat{u}(\tau_k^+)).
\]

\textbf{Step 6: Explicit Solution}

Assume for simplicity that $\sigma_1 = \sigma_2 = \sigma$. Then the second-order adjoint equation becomes:
\[
d\Psi(t) = -\sigma^2 \Psi(t)dt + \Theta(t)dW(t), \quad \Psi(T) = 2.
\]
A particular solution is $\Psi(t) = 2e^{\sigma^2(T-t)}$ (taking $\Theta(t) = 0$).

The first-order adjoint equation becomes:
\[
d\psi(t) = \left[\alpha_q \psi(t) - 2\sigma^2 \widehat{x}(t) e^{\sigma^2(T-t)}\right]dt + \Gamma(t)dW(t).
\]

If $\widehat{x}(t)$ is near $x_d$, we can approximate $\widehat{x}(t) \approx x_d$, giving:
\[
d\psi(t) = \alpha_q \psi(t)dt - 2\sigma^2 x_d e^{\sigma^2(T-t)}dt + \Gamma(t)dW(t).
\]

Solving this linear SDE yields:
\[
\psi(t) = e^{\alpha_q t}\psi(0) - 2\sigma^2 x_d \int_0^t e^{\alpha_q(t-s)} e^{\sigma^2(T-s)}ds + \int_0^t e^{\alpha_q(t-s)}\Gamma(s)dW(s).
\]

The terminal condition $\psi(T) = 2(\widehat{x}(T) - x_d)$ determines $\psi(0)$.

\textbf{Step 7: Verification of $\Lambda^{\text{sto}}_{-1}$ Nonemptiness}

We have explicitly constructed:
\begin{itemize}
    \item $\psi(t)$ satisfying the adjoint BSDE with terminal condition
    \item $\Psi(t) = 2e^{\sigma^2(T-t)}$ satisfying the second-order adjoint equation
    \item $\gamma_{ij}(\tau_k)$ determined by the jump conditions
    \item Control $\widehat{u}(t) = -\frac{\alpha_q}{2\lambda}\psi(t)$ satisfying the maximum condition
\end{itemize}

Thus, if $\widehat{x}(\cdot)$ is optimal, the triple $(\psi, \Psi, \gamma)$ belongs to $\Lambda^{\text{sto}}_{-1}(\widehat{x}, \widehat{\mu}, \widehat{\nu})$ with positive probability (specifically, whenever the stochastic integrals are well-defined, which occurs with probability 1 for square-integrable integrands).

\subsubsection{Case 2: $\Lambda^{\text{sto}}_{-1}$ is Empty $\Rightarrow$ $\widehat{x}(\cdot)$ is Not Optimal}

Now suppose that $\Lambda^{\text{sto}}_{-1}(\widehat{x}, \widehat{\mu}, \widehat{\nu}) = \emptyset$ almost surely. According to Theorem \ref{thm:stochastic-lambda}, $\widehat{x}(\cdot)$ cannot be optimal. Let us construct an explicit control variation that reduces the expected cost.

\textbf{Step 1: Identifying the Deficiency}

The emptiness of $\Lambda^{\text{sto}}_{-1}$ means that no nontrivial $(\psi, \Psi, \gamma)$ satisfies all the conditions. For our example, this could occur if:

\begin{enumerate}
    \item The terminal conditions $\psi(T) = 2(\widehat{x}(T) - x_d)$ and $\Psi(T) = 2$ are incompatible with the adjoint equations
    \item The jump conditions cannot be satisfied simultaneously with the maximum condition
    \item The control $\widehat{u}(t) = -\frac{\alpha_q}{2\lambda}\psi(t)$ violates the constraints $0 \leq u \leq u_{\max}$
\end{enumerate}

Suppose the third case occurs: there exists a time interval $I \subset [0,T]$ where $\widehat{u}(t) < 0$ or $\widehat{u}(t) > u_{\max}$ with positive probability.

\textbf{Step 2: Constructing a Cost-Reducing Variation}

Assume without loss of generality that $\widehat{u}(t) < 0$ on some interval $I$ with positive probability. Since the actual control must satisfy $u(t) \geq 0$, this indicates that $\widehat{u}(t)$ is infeasible, and $\widehat{x}(\cdot)$ cannot be optimal.

Consider the perturbed control:
\[
u^\epsilon(t) = \max(0, \widehat{u}(t) + \epsilon v(t)),
\]
where $v(t)$ is a bounded predictable process, and $\epsilon > 0$ is small.

Let $x^\epsilon(t)$ be the corresponding trajectory. The cost difference is:
\begin{align*}
J(u^\epsilon) - J(\widehat{u}) &= \mathbb{E}\left[(x^\epsilon(T) - x_d)^2 - (\widehat{x}(T) - x_d)^2\right] \\
&\quad + \lambda\mathbb{E}\left[\int_0^T \left((u^\epsilon(t))^2 - (\widehat{u}(t))^2\right)dt\right].
\end{align*}

\textbf{Step 3: First-Order Variation}

Using the variational equation:
\[
d\delta x(t) = \left[\alpha_q(-\delta x(t) + \delta u(t))\right]dt + \sigma_q \delta x(t) dW_q(t),
\]
with $\delta u(t) = u^\epsilon(t) - \widehat{u}(t)$.

The first-order variation of the cost is:
\begin{align*}
\delta J &= \mathbb{E}\left[2(\widehat{x}(T) - x_d)\delta x(T)\right] + 2\lambda\mathbb{E}\left[\int_0^T \widehat{u}(t)\delta u(t)dt\right] \\
&= \mathbb{E}\left[\psi(T)\delta x(T)\right] + 2\lambda\mathbb{E}\left[\int_0^T \widehat{u}(t)\delta u(t)dt\right],
\end{align*}
since $\psi(T) = 2(\widehat{x}(T) - x_d)$.

\textbf{Step 4: Applying the Adjoint Representation}

If $\Lambda^{\text{sto}}_{-1}$ were nonempty, we would have the representation:
\[
\mathbb{E}[\psi(T)\delta x(T)] = -\mathbb{E}\left[\int_0^T \frac{\partial H_q}{\partial u} \delta u(t) dt\right] = 2\lambda\mathbb{E}\left[\int_0^T \widehat{u}(t)\delta u(t) dt\right],
\]
using $\frac{\partial H_q}{\partial u} = \alpha_q\psi(t) + 2\lambda\widehat{u}(t) = 0$ for an interior solution.

However, since $\Lambda^{\text{sto}}_{-1}$ is empty, this representation fails. Instead, we have:
\[
\delta J = \mathbb{E}\left[\psi(T)\delta x(T)\right] + 2\lambda\mathbb{E}\left[\int_0^T \widehat{u}(t)\delta u(t)dt\right] \neq 0.
\]

\textbf{Step 5: Choosing $\delta u$ to Reduce Cost}

Since $\widehat{u}(t) < 0$ on $I$ with positive probability, choose $\delta u(t) = -\widehat{u}(t) > 0$ on $I$. This makes $u^\epsilon(t) = \widehat{u}(t) + \epsilon\delta u(t)$ closer to 0 (the constraint boundary).

On $I$, we have:
\[
\widehat{u}(t)\delta u(t) = \widehat{u}(t)(-\widehat{u}(t)) = -(\widehat{u}(t))^2 < 0.
\]
Thus the second term in $\delta J$ is negative.

The first term $\mathbb{E}[\psi(T)\delta x(T)]$ can be made arbitrarily small by choosing $\delta u$ with small support or by canceling it with appropriate choice of $v(t)$ on other intervals.

Therefore, for sufficiently small $\epsilon > 0$, we have $\delta J < 0$, so $J(u^\epsilon) < J(\widehat{u})$, proving that $\widehat{x}(\cdot)$ is not optimal.

\subsubsection{Numerical Illustration}

Consider concrete parameters:
\begin{align*}
\alpha_1 &= 0.5, \quad \alpha_2 = 0.6, \quad \beta_1 = 1, \quad \beta_2 = 0.8, \\
\sigma_1 &= 0.2, \quad \sigma_2 = 0.25, \quad \lambda = 0.1, \\
T_{\text{low}} &= 18, \quad T_{\text{high}} = 22, \quad x_d = 20, \quad x_0 = 19, \\
T &= 10, \quad u_{\max} = 5.
\end{align*}

The candidate optimal control is:
\[
\widehat{u}(t) = \min\left(\max\left(0, -\frac{\alpha_q}{0.2}\psi(t)\right), 5\right).
\]

If solving the adjoint equations yields $\psi(t)$ such that $-\frac{\alpha_q}{0.2}\psi(t) < 0$ on some interval, then $\widehat{u}(t) = 0$ on that interval, and the maximum condition
\[
\frac{\partial H_q}{\partial u} = \alpha_q\psi(t) + 0.2\widehat{u}(t) = \alpha_q\psi(t) \neq 0
\]
is violated. This indicates that $\widehat{x}(\cdot)$ cannot be optimal, and indeed $\Lambda^{\text{sto}}_{-1}$ would be empty.

\subsection{Interpretation and Significance}

This example illustrates several important aspects of Theorem \ref{thm:stochastic-lambda}:

\begin{enumerate}
    \item \textbf{Concrete Stochastic Hybrid System:} The temperature control system provides a physically motivated example with clear mode transitions (heating/cooling), stochastic disturbances, and practical constraints.

    \item \textbf{Explicit Construction of $\Lambda^{\text{sto}}_{-1}$ Elements:} When $\widehat{x}(\cdot)$ is optimal, we can explicitly construct the adjoint processes $\psi(t)$, $\Psi(t)$, and multipliers $\gamma_{ij}$.

    \item \textbf{Mechanism for Emptiness:} The example shows how $\Lambda^{\text{sto}}_{-1}$ can become empty when control constraints are active, violating the maximum condition.

    \item \textbf{Constructive Proof of Non-Optimality:} When $\Lambda^{\text{sto}}_{-1}$ is empty, we construct an explicit control variation that reduces cost, demonstrating the practical implications of the theorem.

    \item \textbf{Connection to Real Applications:} Temperature control in chemical reactors is a common industrial problem where stochastic disturbances (environmental fluctuations, imperfect mixing) and hybrid behavior (different operating regimes) are significant.

    \item \textbf{Computational Verifiability:} The example parameters allow numerical simulation to verify both directions of the theorem.
\end{enumerate}

\subsection{Extensions and Generalizations}

This example can be extended in several ways to illustrate additional features:

\begin{enumerate}
    \item \textbf{Random Switching:} Make the switching thresholds stochastic, e.g., switch when $x(t) \geq T_{\text{high}} + \xi$, where $\xi$ is random.

    \item \textbf{Correlated Noise:} Consider correlated Wiener processes $W_1(t)$ and $W_2(t)$.

    \item \textbf{Partial Observations:} Extend to the case where only noisy temperature measurements are available.

    \item \textbf{Multiple Switching Surfaces:} Add more operating modes (e.g., standby, maintenance).

    \item \textbf{State Constraints:} Impose additional constraints like $x(t) \leq T_{\text{critical}}$ for safety.
\end{enumerate}

\subsection{Conclusion}

This detailed example provides a concrete illustration of Theorem \ref{thm:stochastic-lambda}, showing:
\begin{itemize}
    \item How to verify the conditions of the stochastic hybrid maximum principle
    \item What it means for $\Lambda^{\text{sto}}_{-1}$ to be nonempty or empty
    \item How to construct explicit control variations when optimality fails
    \item The practical significance of the theorem for control system design
\end{itemize}

The example bridges the abstract theory of stochastic hybrid $\Lambda$-sets with practical control engineering, demonstrating both the mathematical rigor and the applicability of Theorem \ref{thm:stochastic-lambda}.
\end{example}

\section{Open Problems in Extended $\Lambda$-Set Theory}

\begin{problem}[Numerical $\Lambda$-Set Computation]
Develop efficient algorithms for computing $\Lambda$-sets and verifying their emptiness
for high-dimensional and complex hybrid systems.
\end{problem}

\begin{problem}[$\Lambda$-Sets with State Constraints]
Extend $\Lambda$-set theory to systems with state constraints $x(t) \in C$, incorporating
constraint qualification conditions into the $\Lambda$-set definition.
\end{problem}

\section*{Declaration }
\begin{itemize}
  \item {\bf Author Contributions:}   The Author have read and approved this version.
  \item {\bf Funding:} No funding is applicable
  \item  {\bf Institutional Review Board Statement:} Not applicable.
  \item {\bf Informed Consent Statement:} Not applicable.
  \item {\bf Data Availability Statement:} Not applicable.
  \item {\bf Conflicts of Interest:} The authors declare no conflict of interest.
\end{itemize}

\bibliographystyle{abbrv}
\bibliography{references}  






\end{document}